\documentclass{amsart}
\usepackage{amssymb}

\usepackage{amsmath}
\usepackage{amscd}
\usepackage{thmdefs}

\input{tcilatex}
\theoremstyle{definition}
\theoremstyle{remark}
\numberwithin{equation}{section}

\input tcilatex

\begin{document}
\title[Diagonal Compressed Graph $W^{*}$-Probability]{Diagonal Compressed Graph $W^{*}$-Probability}
\author{Ilwoo Cho}
\address{Univ. of Iowa, Dep. of Math, Iowa City, IA, U. S. A}
\email{ilcho@math.uiowa.edu}
\thanks{}
\thanks{}
\maketitle

\begin{abstract}
In this paper, we will consider the compressed graph $W^{*}$-probability
theory. In [15] and [16], we observed Graph $W^{*}$-probability and the
properties of certain amalgamated random variables in the graph $W^{*}$%
-probability space $\left( W^{*}(G),E\right) $ over the diagonal subalgebra $%
D_{G}.$ By using the projections $L_{v},$ $v\in V(G),$ we will consider the
vertex compressed free probability on $\left( W^{*}(G),E\right) .$ Also, for
the fixed vertices $v_{1},...,v_{N}\in V(G),$ we will consider the diagonal
compressed free probability on $\left( W^{*}(G),E\right) .$ We can show that
the diagonal compressed freeness on $\left( W^{*}(G),E\right) $ is preserved
by the $D_{G}$-valued freeness on $\left( W^{*}(G),E\right) .$
\end{abstract}

\strut

In [16], we constructed the graph $W^{*}$-probability spaces. The graph $%
W^{*}$-probability theory is one of the good example of Speicher's
combinatorial free probability theory with amalgamation. In [16], we
observed how to compute the moment and cumulant of an arbitrary random
variables in the graph $W^{*}$-probability space and the freeness on it with
respect to the given conditional expectation. Also, in [17], we consider
certain special random variables of the graph $W^{*}$-probability space, for
example, semicircular elements, even elements and R-diagonal elements. This
shows that the graph $W^{*}$-probability spaces contain the rich free
probabilistic objects. \bigskip Roughly speaking, graph $W^{*}$-algebras are 
$W^{*}$-topology closed version of free semigroupoid algebras defined and
observed by Kribs and Power in [10].

\strut

Throughout this paper, let $G$ be a countable directed graph and let $%
\mathbb{F}^{+}(G)$ be the free semigroupoid of $G,$ in the sense of Kribs
and Power. i.e., it is a collection of all vertices of the graph $G$ as
units and all admissible finite paths, under the admissibility. As a set,
the free semigroupoid $\mathbb{F}^{+}(G)$ can be decomposed by

\strut

\begin{center}
$\mathbb{F}^{+}(G)=V(G)\cup FP(G),$
\end{center}

\strut

where $V(G)$ is the vertex set of the graph $G$ and $FP(G)$ is the set of
all admissible finite paths. Trivially the edge set $E(G)$ of the graph $G$
is properly contained in $FP(G),$ since all edges of the graph can be
regarded as finite paths with their length $1.$ We define a graph $W^{*}$%
-algebra of $G$ by

\strut

\begin{center}
$W^{*}(G)\overset{def}{=}\overline{%
\mathbb{C}[\{L_{w},L_{w}^{*}:w\in
\mathbb{F}^{+}(G)\}]}^{w},$
\end{center}

\strut

where $L_{w}$ and $L_{w}^{\ast }$ are creation operators and annihilation
operators on the generalized Fock space $H_{G}=l^{2}\left( \mathbb{F}%
^{+}(G)\right) $ induced by the given graph $G,$ respectively. Notice that
the creation operators induced by vertices are projections and the creation
operators induced by finite paths are partial isometries. We can define the $%
W^{\ast }$-subalgebra $D_{G}$ of $W^{\ast }(G),$ which is called the
diagonal subalgebra by

\strut

\begin{center}
$D_{G}\overset{def}{=}\overline{\mathbb{C}[\{L_{v}:v\in V(G)\}]}^{w}.$
\end{center}

\strut

Then each element $a$ in the graph $W^{*}$-algebra $W^{*}(G)$ is expressed by

\strut

\begin{center}
$a=\underset{w\in \mathbb{F}^{+}(G:a),\,u_{w}\in \{1,*\}}{\sum }%
p_{w}L_{w}^{u_{w}},$ \ for $p_{w}\in \mathbb{C},$
\end{center}

\strut

where $\mathbb{F}^{+}(G:a)$ is a support of the element $a$, as a subset of
the free semigroupoid $\mathbb{F}^{+}(G).$ The above expression of the
random variable $a$ is said to be the Fourier expansion of $a.$ Since $%
\mathbb{F}^{+}(G)$ is decomposed by the disjoint subsets $V(G)$ and $FP(G),$
the support $\mathbb{F}^{+}(G:a)$ of $a$ is also decomposed by the following
disjoint subsets,

\strut

\begin{center}
$V(G:a)=\mathbb{F}^{+}(G:a)\cap V(G)$
\end{center}

and

\begin{center}
$FP(G:a)=\mathbb{F}^{+}(G:a)\cap FP(G).$
\end{center}

\strut

Thus the operator $a$ can be re-expressed by

\strut

\begin{center}
$a=\underset{v\in V(G:a)}{\sum }p_{v}L_{v}+\underset{w\in FP(G:a),\,u_{w}\in
\{1,*\}}{\sum }p_{w}L_{w}^{u_{w}}.$
\end{center}

\strut

Notice that if $V(G:a)\neq \emptyset ,$ then $\underset{v\in V(G:a)}{\sum }%
p_{v}L_{v}$ is contained in the diagonal subalgebra $D_{G}.$ Thus we have
the canonical conditional expectation $E:W^{*}(G)\rightarrow D_{G},$ defined
by

\strut

\begin{center}
$E\left( a\right) =\underset{v\in V(G:a)}{\sum }p_{v}L_{v},$
\end{center}

\strut

for all $a=\underset{w\in \mathbb{F}^{+}(G:a),\,u_{w}\in \{1,*\}}{\sum }%
p_{w}L_{w}^{u_{w}}$ \ in $W^{*}(G).$ Then the algebraic pair $\left(
W^{*}(G),E\right) $ is a $W^{*}$-probability space with amalgamation over $%
D_{G}$ (See [16]). It is easy to check that the conditional expectation $E$
is faithful in the sense that if $E(a^{*}a)=0_{D_{G}},$ for $a\in W^{*}(G),$
then $a=0_{D_{G}}.$

\strut

For the fixed operator $a\in W^{*}(G),$ the support $\mathbb{F}^{+}(G:a)$ of
the operator $a$ is again decomposed by

\strut

\begin{center}
$\mathbb{F}^{+}(G:a)=V(G:a)\cup FP_{*}(G:a)\cup FP_{*}^{c}(G:a),$
\end{center}

\strut

with the decomposition of $FP(G:a),$

$\strut $

\begin{center}
$FP(G:a)=FP_{*}(G:a)\cup FP_{*}^{c}(G:a),$
\end{center}

where

\strut

\begin{center}
$FP_{*}(G:a)=\{w\in FP(G:a):$both $L_{w}$ and $L_{w}^{*}$ are summands of $%
a\}$
\end{center}

and

\begin{center}
$FP_{*}(G:a)=FP(G:a)\,\,\setminus \,\,FP_{*}(G:a).$
\end{center}

\strut

The above new expression plays a key role to find the $D_{G}$-valued moments
of the random variable $a.$ In fact, the summands $p_{v}L_{v}$'s and $%
p_{w}L_{w}+p_{w^{t}}L_{w}^{*},$ for $v\in V(G:a)$ and $w\in FP_{*}(G:a)$ act
for the computation of $D_{G}$-valued moments of $a.$ By using the above
partition of the support of a random variable, we can compute the $D_{G}$%
-valued moments and $D_{G}$-valued cumulants of it via the lattice path
model $LP_{n}$ and the lattice path model $LP_{n}^{*}$ satisfying the $*$%
-axis-property. At a first glance, the computations of $D_{G}$-valued
moments and cumulants look so abstract and hence it looks useless. However,
these computations, in particular the computation of $D_{G}$-valued
cumulants, provides us how to figure out the $D_{G}$-freeness of random
variables by making us compute the mixed cumulants. As applications, in the
final chapter, we can compute the moment and cumulant of the operator that
is the sum of $N$-free semicircular elements with their covariance $2.$

\strut \strut \strut

Based on the $D_{G}$-cumulant computation, we can characterize the $D_{G}$%
-freeness of generators of $W^{*}(G),$ by the so-called diagram-distinctness
on the graph $G.$ i.e., the random variables $L_{w_{1}}$ and $L_{w_{2}}$ are
free over $D_{G}$ if and only if $w_{1}$ and $w_{2}$ are diagram-distinct
the sense that $w_{1}$ and $w_{2}$ have different diagrams on the graph $G.$
Also, we could find the necessary condition for the $D_{G}$-freeness of two
arbitrary random variables $a$ and $b.$ i.e., if the supports $\mathbb{F}%
^{+}(G:a)$ and $\mathbb{F}^{+}(G:b)$ are diagram-distinct, in the sense that 
$w_{1}$ and $w_{2}$ are diagram distinct for all pairs $(w_{1},w_{2})$ $\in $
$\mathbb{F}^{+}(G:a)$ $\times $ $\mathbb{F}^{+}(G:b),$ then the random
variables $a$ and $b$ are free over $D_{G}.$

\strut \strut

In [17], we considered some special $D_{G}$-valued random variables in a
graph $W^{*}$-probability space $\left( W^{*}(G),E\right) .$ The those
random variables are the basic objects to study Free Probability Theory. We
can conclude that

\strut

(i) \ \ if $l$ is a loop, then $L_{l}+L_{l}^{*}$ is $D_{G}$-semicircular.

\strut

(ii) \ if $w$ is a finite path, then $L_{w}+L_{w}^{*}$ is $D_{G}$-even.

\strut

(iii) if $w$ is a finite path, then $L_{w}$ and $L_{w}^{*}$ are $D_{G}$%
-valued R-diagonal.

\strut

In this paper, we will observe the diagonal compressed random variables in
the graph $W^{*}$-probability space $\left( W^{*}(G),E\right) .$ Let $%
v_{1},...,v_{N}\in V(G)$ and let $a$ be a $D_{G}$-valued random variable in $%
\left( W^{*}(G),E\right) .$ Define the diagonal compressed random variable
of $a$ by $V=\{v_{1},...,v_{N}\}$ by

\strut

\begin{center}
$L_{v_{1}}aL_{v_{1}}+...+L_{v_{N}}aL_{v_{N}}.$
\end{center}

\strut

Notice that if $v\in V(G),$ then $L_{v}aL_{v}$ is the compressed random
variable by $L_{v}$ and the compressed random variable has its support
contained in $\{v\}$ $\cup $ $loop_{v}(G),$ where $loop_{v}(G)$ $=$ $\{l\in
loop(G):$ $l=vlv\}.$ We will consider the $D_{G}$-moments, $D_{G}$-cumulants
and $D_{G}$-freeness of such compressed random variables.

\strut \strut

\strut \strut

\strut \strut

\section{Graph $W^{*}$-Probability Theory}

\strut

\strut

Let $G$ be a countable directed graph and let $\Bbb{F}^{+}(G)$ be the free
semigroupoid of $G.$ i.e., the set $\mathbb{F}^{+}(G)$ is the collection of
all vertices as units and all admissible finite paths of $G.$ Let $w$ be a
finite path with its source $s(w)=x$ and its range $r(w)=y,$ where $x,y\in
V(G).$ Then sometimes we will denote $w$ by $w=xwy$ to express the source
and the range of $w.$ We can define the graph Hilbert space $H_{G}$ by the
Hilbert space $l^{2}\left( \mathbb{F}^{+}(G)\right) $ generated by the
elements in the free semigroupoid $\mathbb{F}^{+}(G).$ i.e., this Hilbert
space has its Hilbert basis $\mathcal{B}=\{\xi _{w}:w\in \mathbb{F}%
^{+}(G)\}. $ Suppose that $w=e_{1}...e_{k}\in FP(G)$ is a finite path with $%
e_{1},...,e_{k}\in E(G).$ Then we can regard $\xi _{w}$ as $\xi
_{e_{1}}\otimes ...\otimes \xi _{e_{k}}.$ So, in [10], Kribs and Power
called this graph Hilbert space the generalized Fock space. Throughout this
paper, we will call $H_{G}$ the graph Hilbert space to emphasize that this
Hilbert space is induced by the graph.

\strut

Define the creation operator $L_{w},$ for $w\in \mathbb{F}^{+}(G),$ by the
multiplication operator by $\xi _{w}$ on $H_{G}.$ Then the creation operator 
$L$ on $H_{G}$ satisfies that

\strut

(i) \ $L_{w}=L_{xwy}=L_{x}L_{w}L_{y},$ for $w=xwy$ with $x,y\in V(G).$

\strut

(ii) $L_{w_{1}}L_{w_{2}}=\left\{ 
\begin{array}{lll}
L_{w_{1}w_{2}} &  & \text{if }w_{1}w_{2}\in \mathbb{F}^{+}(G) \\ 
&  &  \\ 
0 &  & \text{if }w_{1}w_{2}\notin \mathbb{F}^{+}(G),
\end{array}
\right. $

\strut

\ \ \ \ for all $w_{1},w_{2}\in \mathbb{F}^{+}(G).$

\strut

Now, define the annihilation operator $L_{w}^{*},$ for $w\in \mathbb{F}%
^{+}(G)$ by

\strut

\begin{center}
$L_{w}^{\ast }\xi _{w^{\prime }}\overset{def}{=}\left\{ 
\begin{array}{lll}
\xi _{h} &  & \text{if }w^{\prime }=wh\in \mathbb{F}^{+}(G)\xi \\ 
&  &  \\ 
0 &  & \text{otherwise.}
\end{array}
\right. $
\end{center}

\strut

The above definition is gotten by the following observation ;

\strut

\begin{center}
$
\begin{array}{ll}
<L_{w}\xi _{h},\xi _{wh}>\, & =\,<\xi _{wh},\xi _{wh}>\, \\ 
& =\,1=\,<\xi _{h},\xi _{h}> \\ 
& =\,<\xi _{h},L_{w}^{*}\xi _{wh}>,
\end{array}
\,$
\end{center}

\strut

where $<,>$ is the inner product on the graph Hilbert space $H_{G}.$ Of
course, in the above formula we need the admissibility of $w$ and $h$ in $%
\mathbb{F}^{+}(G).$ However, even though $w$ and $h$ are not admissible
(i.e., $wh\notin \mathbb{F}^{+}(G)$), by the definition of $L_{w}^{\ast },$
we have that

\strut

\begin{center}
$
\begin{array}{ll}
<L_{w}\xi _{h},\xi _{h}> & =\,<0,\xi _{h}> \\ 
& =0=\,<\xi _{h},0> \\ 
& =\,<\xi _{h},L_{w}^{*}\xi _{h}>.
\end{array}
\,\,$
\end{center}

\strut

Notice that the creation operator $L$ and the annihilation operator $L^{*}$
satisfy that

\strut

(1.1) \ \ \ $L_{w}^{*}L_{w}=L_{y}$ \ \ and \ \ $L_{w}L_{w}^{*}=L_{x},$ \ for
all \ $w=xwy\in \mathbb{F}^{+}(G),$

\strut

\textbf{under the weak topology}, where $x,y\in V(G).$ Remark that if we
consider the von Neumann algebra $W^{*}(\{L_{w}\})$ generated by $L_{w}$ and 
$L_{w}^{*}$ in $B(H_{G}),$ then the projections $L_{y}$ and $L_{x}$ are
Murray-von Neumann equivalent, because there exists a partial isometry $%
L_{w} $ satisfying the relation (1.1). Indeed, if $w=xwy$ in $\mathbb{F}%
^{+}(G), $ with $x,y\in V(G),$ then under the weak topology we have that

\strut

(1,2) \ \ \ $L_{w}L_{w}^{*}L_{w}=L_{w}$ \ \ and \ \ $%
L_{w}^{*}L_{w}L_{w}^{*}=L_{w}^{*}.$

\strut

So, the creation operator $L_{w}$ is a partial isometry in $W^{*}(\{L_{w}\})$
in $B(H_{G}).$ Assume now that $v\in V(G).$ Then we can regard $v$ as $%
v=vvv. $ So,

\strut

(1.3) $\ \ \ \ \ \ \ \ \ L_{v}^{*}L_{v}=L_{v}=L_{v}L_{v}^{*}=L_{v}^{*}.$

\strut

This relation shows that $L_{v}$ is a projection in $B(H_{G})$ for all $v\in
V(G).$

\strut

Define the \textbf{graph }$W^{*}$\textbf{-algebra} $W^{*}(G)$ by

\strut

\begin{center}
$W^{*}(G)\overset{def}{=}\overline{%
\mathbb{C}[\{L_{w},L_{w}^{*}:w\in
\mathbb{F}^{+}(G)\}]}^{w}.$
\end{center}

\strut

Then all generators are either partial isometries or projections, by (1.2)
and (1.3). So, this graph $W^{\ast }$-algebra contains a rich structure, as
a von Neumann algebra. (This construction can be the generalization of that
of group von Neumann algebra.) Naturally, we can define a von Neumann
subalgebra $D_{G}\subset W^{\ast }(G)$ generated by all projections $L_{v},$ 
$v\in V(G).$ i.e.

\strut

\begin{center}
$D_{G}\overset{def}{=}W^{*}\left( \{L_{v}:v\in V(G)\}\right) .$
\end{center}

\strut

We call this subalgebra the \textbf{diagonal subalgebra} of $W^{*}(G).$
Notice that $D_{G}=\Delta _{\left| G\right| }\subset M_{\left| G\right| }(%
\mathbb{C}),$ where $\Delta _{\left| G\right| }$ is the subalgebra of $%
M_{\left| G\right| }(\mathbb{C})$ generated by all diagonal matrices. Also,
notice that $1_{D_{G}}=\underset{v\in V(G)}{\sum }L_{v}=1_{W^{*}(G)}.$

\strut

If $a\in W^{*}(G)$ is an operator, then it has the following decomposition
which is called the Fourier expansion of $a$ ;

\strut

(1.4) $\ \ \ \ \ \ \ \ \ \ \ a=\underset{w\in \mathbb{F}^{+}(G:a),\,u_{w}\in
\{1,*\}}{\sum }p_{w}L_{w}^{u_{w}},$

\strut

where $p_{w}\in C$ and $\mathbb{F}^{+}(G:a)$ is the support of $a$ defined by

\strut

\begin{center}
$\mathbb{F}^{+}(G:a)=\{w\in \mathbb{F}^{+}(G):p_{w}\neq 0\}.$
\end{center}

\strut

Remark that the free semigroupoid $\mathbb{F}^{+}(G)$ has its partition $%
\{V(G),FP(G)\},$ as a set. i.e.,

\strut

\begin{center}
$\mathbb{F}^{+}(G)=V(G)\cup FP(G)$ \ \ and \ \ $V(G)\cap FP(G)=\emptyset .$
\end{center}

\strut

So, the support of $a$ is also partitioned by

\strut

\begin{center}
$\mathbb{F}^{+}(G:a)=V(G:a)\cup FP(G:a),$
\end{center}

\strut where

\begin{center}
$V(G:a)\overset{def}{=}V(G)\cap \mathbb{F}^{+}(G:a)$
\end{center}

and

\begin{center}
$FP(G:a)\overset{def}{=}FP(G)\cap \mathbb{F}^{+}(G:a).$
\end{center}

\strut

So, the above Fourier expansion (1.4) of the random variable $a$ can be
re-expressed by

\strut

(1.5) $\ \ \ \ \ \ a=\underset{v\in V(G:a)}{\sum }p_{v}L_{v}+\underset{w\in
FP(G:a),\,u_{w}\in \{1,*\}}{\sum }p_{w}L_{w}^{u_{w}}.$

\strut

We can easily see that if $V(G:a)\neq \emptyset ,$ then $\underset{v\in
V(G:a)}{\sum }p_{v}L_{v}$ is contained in the diagonal subalgebra $D_{G}.$
Also, if $V(G:a)=\emptyset ,$ then $\underset{v\in V(G:a)}{\sum }%
p_{v}L_{v}=0_{D_{G}}.$ So, we can define the following canonical conditional
expectation $E:W^{*}(G)\rightarrow D_{G}$ by

\strut

(1.6) \ \ \ $E(a)=E\left( \underset{w\in \mathbb{F}^{+}(G:a),\,u_{w}\in
\{1,*\}}{\sum }p_{w}L_{w}^{u_{w}}\right) \overset{def}{=}\underset{v\in
V(G:a)}{\sum }p_{v}L_{v},$

\strut

for all $a\in W^{*}(G).$ Indeed, $E$ is a well-determined conditional
expectation.

\strut \strut \strut \strut

\begin{definition}
Let $G$ be a countable directed graph and let $W^{*}(G)$ be the graph $W^{*}$%
-algebra induced by $G.$ Let $E:W^{*}(G)\rightarrow D_{G}$ be the
conditional expectation defined above. Then we say that the algebraic pair $%
\left( W^{*}(G),E\right) $ is the graph $W^{*}$-probability space over the
diagonal subalgebra $D_{G}$. By the very definition, it is one of the $W^{*}$%
-probability space with amalgamation over $D_{G}.$ All elements in $\left(
W^{*}(G),E\right) $ are called $D_{G}$-valued random variables.
\end{definition}

\strut

We have a graph $W^{*}$-probability space $\left( W^{*}(G),E\right) $ over
its diagonal subalgebra $D_{G}.$ We will define the following free
probability data of $D_{G}$-valued random variables.

\strut

\begin{definition}
Let $W^{*}(G)$ be the graph $W^{*}$-algebra induced by $G$ and let $a\in
W^{*}(G).$ Define the $n$-th ($D_{G}$-valued) moment of $a$ by

\strut 

$\ \ \ \ \ E\left( d_{1}ad_{2}a...d_{n}a\right) ,$ for all $n\in \mathbb{N}$,

\strut 

where $d_{1},...,d_{n}\in D_{G}$. Also, define the $n$-th ($D_{G}$-valued)
cumulant of $a$ by

\strut 

$\ \ \ \ \ k_{n}(d_{1}a,d_{2}a,...,d_{n}a)=C^{(n)}\left( d_{1}a\otimes
d_{2}a\otimes ...\otimes d_{n}a\right) ,$

\strut 

for all $n\in \mathbb{N},$ and for $d_{1},...,d_{n}\in D_{G},$ where $%
\widehat{C}=(C^{(n)})_{n=1}^{\infty }\in I^{c}\left( W^{*}(G),D_{G}\right) $
is the cumulant multiplicative bimodule map induced by the conditional
expectation $E,$ in the sense of Speicher. We define the $n$-th trivial
moment of $a$ and the $n$-th trivial cumulant of $a$ by

\strut 

$\ \ \ \ \ E(a^{n})$ $\ \ $and $\ \ k_{n}\left( \underset{n-times}{%
\underbrace{a,a,...,a}}\right) =C^{(n)}\left( a\otimes a\otimes ...\otimes
a\right) ,$

\strut 

respectively, for all $n\in \mathbb{N}.$
\end{definition}

\strut

To compute the $D_{G}$-valued moments and cumulants of the $D_{G}$-valued
random variable $a,$ we need to introduce the following new definition ;

\strut

\begin{definition}
Let $\left( W^{*}(G),E\right) $ be a graph $W^{*}$-probability space over $%
D_{G}$ and let $a\in \left( W^{*}(G),E\right) $ be a random variable. Define
the subset $FP_{*}(G:a)$ in $FP(G:a)$ \ by

\strut 

$\ \ \ FP_{*}\left( G:a\right) \overset{def}{=}\{w\in \mathbb{F}^{+}(G:a):$%
both $L_{w}$ and $L_{w}^{*}$ are summands of $a\}.$

\strut 

And let $FP_{*}^{c}(G:a)\overset{def}{=}FP(G:a)\,\setminus \,FP_{*}(G:a).$
\end{definition}

\strut \strut \strut

We already observed that if $a\in \left( W^{*}(G),E\right) $ is a $D_{G}$%
-valued random variable, then $a$ has its Fourier expansion $a_{d}+a_{0},$
where

\strut

\begin{center}
$a_{d}=\underset{v\in V(G:a)}{\sum }p_{v}L_{v}$
\end{center}

and

\begin{center}
$a_{0}=\underset{w\in FP(G:a),\,u_{w}\in \{1,*\}}{\sum }p_{w}L_{w}^{u_{w}}.$
\end{center}

\strut

By the previous definition, the set $FP(G:a)$ is partitioned by

\strut

\begin{center}
$FP(G:a)=FP_{*}(G:a)\cup FP_{*}^{c}(G:a),$
\end{center}

\strut

for the fixed random variable $a$ in $\left( W^{*}(G),E\right) .$ So, the
summand $a_{0},$ in the Fourier expansion of $a=a_{d}+a_{0},$ has the
following decomposition ;

\strut

\begin{center}
$a_{0}=a_{(*)}+a_{(non-*)},$
\end{center}

\strut where\strut

\begin{center}
$a_{(*)}=\underset{l\in FP_{*}(G:a)}{\sum }\left(
p_{l}L_{l}+p_{l^{t}}L_{l}^{*}\right) $
\end{center}

and

\begin{center}
$a_{(non-*)}=\underset{w\in FP_{*}^{c}(G:a),\,u_{w}\in \{1,*\}}{\sum }%
p_{w}L_{w}^{u_{w}},$
\end{center}

\strut

where $p_{l^{t}}$ is the coefficient of $L_{l}^{*}$ depending on $l\in
FP_{*}(G:a).$

\strut \strut \strut

\strut \strut

\strut \strut

\subsection{$D_{G}$-Moments and $D_{G}$-Cumulants of Random Variables}

\strut

\strut

\strut

Throughout this chapter, let $G$ be a countable directed graph and let $%
\left( W^{*}(G),E\right) $ be the graph $W^{*}$-probability space over its
diagonal subalgebra $D_{G}.$ In this chapter, we will compute the $D_{G}$%
-valued moments and the $D_{G}$-valued cumulants of arbitrary random variable

$\strut $

\begin{center}
$a=\underset{w\in \mathbb{F}^{+}(G:a),\,u_{w}\in \{1,*\}}{\sum }%
p_{w}L_{w}^{u_{w}}$
\end{center}

\strut

in the graph $W^{*}$-probability space $\left( W^{*}(G),E\right) $.

\strut

\strut

\subsubsection{Lattice Path Model}

\strut

\strut

\strut

Throughout this section, let $G$ be a countable directed graph and let $%
\left( W^{*}(G),E\right) $ be the graph $W^{*}$-probability space over its
diagonal subalgebra $D_{G}.$ Let $w_{1},...,w_{n}\in \Bbb{F}^{+}(G)$ and let 
$L_{w_{1}}^{u_{w_{1}}}...L_{w_{n}}^{u_{w_{n}}}\in \left( W^{*}(G),E\right) $
be a $D_{G}$-valued random variable. In this section, we will define a
lattice path model for the random variable $%
L_{w_{1}}^{u_{w_{1}}}...L_{w_{n}}^{u_{w_{n}}}.$ Recall that if $%
w=e_{1}....e_{k}\in FP(G)$ with $e_{1},...,e_{k}\in E(G),$ then we can
define the length $\left| w\right| $ of $w$ by $k.$ i.e.e, the length $%
\left| w\right| $ of $w$ is the cardinality $k$ of the admissible edges $%
e_{1},...,e_{k}.$

\strut

\begin{definition}
Let $G$ be a countable directed graph and $\Bbb{F}^{+}(G),$ the free
semigroupoid. If $w\in \Bbb{F}^{+}(G),$ then $L_{w}$ is the corresponding $%
D_{G}$-valued random variable in $\left( W^{*}(G),E\right) .$ We define the
lattice path $l_{w}$ of $L_{w}$ and the lattice path $l_{w}^{-1}$ of $%
L_{w}^{*}$ by the lattice paths satisfying that ;

\strut 

(i) \ \ the lattice path $l_{w}$ starts from $*=(0,0)$ on the $\Bbb{R}^{2}$%
-plane.

\strut 

(ii) \ if $w\in V(G),$ then $l_{w}$ has its end point $(0,1).$

\strut 

(iii) if $w\in E(G),$ then $l_{w}$ has its end point $(1,1).$

\strut 

(iv) if $w\in E(G),$ then $l_{w}^{-1}$ has its end point $(-1,-1).$

\strut 

(v) \ if $w\in FP(G)$ with $\left| w\right| =k,$ then $l_{w}$ has its end
point $(k,k).$

\strut 

(vi) if $w\in FP(G)$ with $\left| w\right| =k,$ then $l_{w}^{-1}$ has its
end point $(-k,-k).$

\strut 

Assume that finite paths $w_{1},...,w_{s}$ in $FP(G)$ satisfy that $%
w_{1}...w_{s}\in FP(G).$ Define the lattice path $l_{w_{1}...w_{s}}$ by the
connected lattice path of the lattice paths $l_{w_{1}},$ ..., $l_{w_{s}}.$
i.e.e, $l_{w_{2}}$ starts from $(k_{w_{1}},k_{w_{1}})\in \Bbb{R}^{+}$ and
ends at $(k_{w_{1}}+k_{w_{2}},k_{w_{1}}+k_{w_{2}}),$ where $\left|
w_{1}\right| =k_{w_{1}}$ and $\left| w_{2}\right| =k_{w_{2}}.$ Similarly, we
can define the lattice path $l_{w_{1}...w_{s}}^{-1}$ as the connected path
of $l_{w_{s}}^{-1},$ $l_{w_{s-1}}^{-1},$ ..., $l_{w_{1}}^{-1}.$
\end{definition}

\strut

\begin{definition}
Let $G$ be a countable directed graph and assume that $%
L_{w_{1}},...,L_{w_{n}}$ are generators of $\left( W^{*}(G),E\right) .$ Then
we have the lattice paths $l_{w_{1}},$ ..., $l_{w_{n}}$ of $L_{w_{1}},$ ..., 
$L_{w_{n}},$ respectively in $\Bbb{R}^{2}.$ Suppose that $%
L_{w_{1}}^{u_{w_{1}}}...L_{w_{n}}^{u_{w_{n}}}\neq 0_{D_{G}}$ in $\left(
W^{*}(G),E\right) ,$ where $u_{w_{1}},...,u_{w_{n}}\in \{1,*\}.$ Define the
lattice path $l_{w_{1},...,w_{n}}^{u_{w_{1}},...,u_{w_{n}}}$ of nonzero $%
L_{w_{1}}^{u_{w_{1}}}...L_{w_{n}}^{u_{w_{n}}}$ by the connected lattice path
of $l_{w_{1}}^{t_{w_{1}}},$ ..., $l_{w_{n}}^{t_{w_{n}}},$ where $t_{w_{j}}=1$
if $u_{w_{j}}=1$ and $t_{w_{j}}=-1$ if $u_{w_{j}}=*.$ Assume that $%
L_{w_{1}}^{u_{w_{1}}}...L_{w_{n}}^{u_{w_{n}}}$ $=$ $0_{D_{G}}.$ Then the
empty set $\emptyset $ in $\Bbb{R}^{2}$ is the lattice path of it. We call
it the empty lattice path. By $LP_{n},$ we will denote the set of all
lattice paths of the $D_{G}$-valued random variables having their forms of $%
L_{w_{1}}^{u_{w_{1}}}...L_{w_{n}}^{u_{w_{n}}},$ including empty lattice path.
\end{definition}

\strut

Also, we will define the following important property on the set of all
lattice paths ;

\strut

\begin{definition}
Let $l_{w_{1},...,w_{n}}^{u_{w_{1}},...,u_{w_{n}}}\neq \emptyset $ be a
lattice path of $L_{w_{1}}^{u_{w_{1}}}...L_{w_{n}}^{u_{w_{n}}}\neq 0_{D_{G}}$
in $LP_{n}.$ If the lattice path $%
l_{w_{1},...,w_{n}}^{u_{w_{1}},...,u_{w_{n}}}$ starts from $*$ and ends on
the $*$-axis in $\Bbb{R}^{+},$ then we say that the lattice path $%
l_{w_{1},...,w_{n}}^{u_{w_{1}},...,u_{w_{n}}}$ has the $*$-axis-property. By 
$LP_{n}^{*},$ we will denote the set of all lattice paths having their forms
of $l_{w_{1},...,w_{n}}^{u_{w_{1}},...,u_{w_{n}}}$ which have the $*$%
-axis-property. By little abuse of notation, sometimes, we will say that the 
$D_{G}$-valued random variable $L_{w_{1}}^{u_{w_{1}}}...L_{w_{n}}^{u_{w_{n}}}
$satisfies the $*$-axis-property if the lattice path $%
l_{w_{1},...,w_{n}}^{u_{w_{1}},...,u_{w_{n}}}$ of it has the $*$%
-axis-property.
\end{definition}

\strut

The following theorem shows that finding $E\left(
L_{w_{1}}^{u_{w_{1}}}...L_{w_{n}}^{u_{w_{n}}}\right) $ is checking the $*$%
-axis-property of $L_{w_{1}}^{u_{w_{1}}}...L_{w_{n}}^{u_{w_{n}}}.$

\strut

\begin{theorem}
(See [15]) Let $L_{w_{1}}^{u_{w_{1}}}...L_{w_{n}}^{u_{w_{n}}}\in \left(
W^{*}(G),E\right) $ be a $D_{G}$-valued random variable, where $%
u_{w_{1}},...,u_{w_{n}}\in \{1,*\}.$ Then $E\left(
L_{w_{1}}^{u_{w_{1}}}...L_{w_{n}}^{u_{w_{n}}}\right) $ $\neq $ $0_{D_{G}}$
if and only if $L_{w_{1}}^{u_{w_{1}}}...L_{w_{n}}^{u_{w_{n}}}$ has the $*$%
-axis-property (i.e., the corresponding lattice path $%
l_{w_{1},...,w_{n}}^{u_{w_{1}},...,u_{w_{n}}}$ of $%
L_{w_{1}}^{u_{w_{1}}}...L_{w_{n}}^{u_{w_{n}}}$ is contained in $LP_{n}^{*}.$
Notice that $\emptyset \notin LP_{n}^{*}.$) \ $\square $
\end{theorem}

\strut \strut 

By the previous theorem, we can conclude that $E\left(
L_{w_{1}}^{u_{w_{1}}}...L_{w_{n}}^{u_{w_{n}}}\right) =L_{v},$ for some $v\in
V(G)$ if and only if the lattice path $%
l_{w_{1},...,w_{n}}^{u_{w_{1}},...,u_{w_{n}}}$ has the $*$-axis-property
(i.e., $l_{w_{1},...,w_{n}}^{u_{w_{1}},...,u_{w_{n}}}\in LP_{n}^{*}$).\strut
\strut \strut

\strut 

\strut \strut

\subsubsection{$D_{G}$-Valued Moments and Cumulants of Random Variables\strut
}

\bigskip

\bigskip

Let $w_{1},...,w_{n}\in \mathbb{F}^{+}(G)$, $u_{1},...,u_{n}\in \{1,*\}$ and
let $L_{w_{1}}^{u_{1}}...L_{w_{n}}^{u_{n}}\in \left( W^{*}(G),E\right) $ be
a $D_{G}$-valued random variable. Recall that, in the previous section, we
observed that the $D_{G}$-valued random variable $%
L_{w_{1}}^{u_{1}}...L_{w_{n}}^{u_{n}}=L_{v}\in \left( W^{*}(G),E\right) $
with $v\in V(G)$ if and only if the lattice path $%
l_{w_{1},...,w_{n}}^{u_{1},...,u_{n}}$ of \ $%
L_{w_{1}}^{u_{1}}...L_{w_{n}}^{u_{n}}$ has the $*$-axis-property
(equivalently, $l_{w_{1},...,w_{n}}^{u_{1},...,u_{n}}\in LP_{n}^{*}$).
Throughout this section, fix a $D_{G}$-valued random variable $a\in \left(
W^{*}(G),E\right) .$ Then the $D_{G}$-valued random variable $a$ has the
following Fourier expansion,

\bigskip

\begin{center}
$a=\underset{v\in V(G:a)}{\sum }p_{v}L_{v}+\underset{l\in FP_{*}(G:a)}{\sum }%
\left( p_{l}L_{l}+p_{l^{t}}L_{l}\right) +\underset{w\in
FP_{*}^{c}(G:a),~u_{w}\in \{1,*\}}{\sum }p_{w}L_{w}^{u_{w}}.$
\end{center}

\bigskip

Let's observe the new $D_{G}$-valued random variable $d_{1}ad_{2}a...d_{n}a%
\in \left( W^{*}(G),E\right) ,$ where $d_{1},...,d_{n}\in D_{G}$ and $a\in
W^{*}(G)$ is given. Put

\strut

\begin{center}
$d_{j}=\underset{v_{j}\in V(G:d_{j})}{\sum }q_{v_{j}}L_{v_{j}}\in D_{G},$
for \ $j=1,...,n.$
\end{center}

\strut

Notice that $V(G:d_{j})=\mathbb{F}^{+}(G:d_{j}),$ since $d_{j}\in
D_{G}\hookrightarrow W^{\ast }(G).$ Then

\strut

$\ d_{1}ad_{2}a...d_{n}a$

\strut

$\ \ \ =\left( \underset{v_{1}\in V(G:d_{1})}{\sum }q_{v_{1}}L_{v_{1}}%
\right) \left( \underset{w_{1}\in \mathbb{F}^{+}(G:a),\,u_{w_{1}}\in \{1,*\}%
}{\sum }p_{w_{1}}L_{w_{1}}^{u_{w_{1}}}\right) $

$\ \ \ \ \ \ \ \ \ \ \ \ \ \ \ \ \ \cdot \cdot \cdot \left( \underset{%
v_{1}\in V(G:d_{n})}{\sum }q_{v_{n}}L_{v_{n}}\right) \left( \underset{%
w_{n}\in \mathbb{F}^{+}(G:a),\,u_{w_{n}}\in \{1,*\}}{\sum }%
p_{w_{n}}L_{w_{n}}^{u_{w_{n}}}\right) $

\strut

$\ \ \ =\underset{(v_{1},...,v_{n})\in \Pi _{j=1}^{n}V(G:d_{j})}{\sum }%
\left( q_{v_{1}}...q_{v_{n}}\right) $

$\ \ \ \ \ \ \ \ \ \ \ \ \ \ \ \ \ \ \ \ \ \ \ \ \ \ \ \ \ \ \ \
(L_{v_{1}}\left( \underset{w_{1}\in \mathbb{F}^{+}(G:a),\,u_{w_{1}}\in
\{1,*\}}{\sum }p_{w_{1}}L_{w_{1}}^{u_{w_{1}}}\right) $

$\ \ \ \ \ \ \ \ \ \ \ \ \ \ \ \ \ \ \ \ \ \ \ \ \ \ \ \ \ \ \ \ \ \ \ \ \ \
\ \cdot \cdot \cdot L_{v_{n}}\left( \underset{w_{n}\in \mathbb{F}%
^{+}(G:a),\,u_{w_{n}}\in \{1,*\}}{\sum }p_{w_{n}}L_{w_{n}}^{u_{w_{n}}}%
\right) )$

\strut

(1.7)

\strut

$\ \ \ =\underset{(v_{1},...,v_{n})\in \Pi _{j=1}^{n}V(G:d_{j})}{\sum }%
\left( q_{v_{1}}...q_{v_{n}}\right) $

\strut

$\ \ \ \ \ \underset{(w_{1},...,w_{n})\in \mathbb{F}^{+}(G:a)^{n},\,%
\,u_{w_{j}}\in \{1,*\}}{\sum }\left( p_{w_{1}}...p_{w_{n}}\right)
L_{v_{1}}L_{w_{1}}^{u_{w_{1}}}...L_{v_{n}}L_{w_{n}}^{u_{w_{n}}}.$

\strut

Now, consider the random variable $%
L_{v_{1}}L_{w_{1}}^{u_{w_{1}}}...L_{v_{n}}L_{w_{n}}^{u_{w_{n}}}$ in the
formula (1.11). Suppose that $w_{j}=x_{j}w_{j}y_{j},$ with $x_{j},y_{j}\in
V(G),$ for all \ $j=1,...,n.$ Then

\strut

\strut (1.8)

\begin{center}
$L_{v_{1}}L_{w_{1}}^{u_{w_{1}}}...L_{v_{n}}L_{w_{n}}^{u_{w_{n}}}=\delta
_{(v_{1},x_{1},y_{1}:u_{w_{1}})}\cdot \cdot \cdot \delta
_{(v_{n},x_{n},y_{n}:u_{w_{n}})}\left(
L_{w_{1}}^{u_{w_{n}}}...L_{w_{n}}^{u_{w_{n}}}\right) ,$
\end{center}

\strut

where

\strut

\begin{center}
$\delta _{(v_{j},x_{j},y_{j}:u_{w_{j}})}=\left\{ 
\begin{array}{lll}
\delta _{v_{j},x_{j}} &  & \text{if }u_{w_{j}}=1 \\ 
&  &  \\ 
\delta _{v_{j},y_{j}} &  & \text{if }u_{w_{j}}=*,
\end{array}
\right. $
\end{center}

\strut

for all \ $j=1,...,n,$ where $\delta $ in the right-hand side is the
Kronecker delta. So, the left-hand side can be understood as a (conditional)
Kronecker delta depending on $\{1,*\}$.

\strut \strut

By (1.7) and (1.8), the $n$-th moment of $a$ is determined by ;

\strut \strut 

\begin{proposition}
Let $a\in \left( W^{*}(G),E\right) $ be given as above. Then the $n$-th
moment of $a$ is

\strut 

$\ \ E\left( d_{1}a...d_{n}a\right) =\underset{(v_{1},...,v_{n})\in \Pi
_{j=1}^{n}V(G:d_{j})}{\sum }\left( \Pi _{j=1}^{n}q_{v_{j}}\right) $

\strut 

$\ \ \ \underset{(w_{1},...,w_{n})\in \mathbb{F}^{+}(G:a)^{n},\,u_{w_{j}}\in
\{1,*\},\,l_{w_{1},...,w_{n}}^{u_{w_{1}},...,u_{w_{n}}}\in LP_{n}^{*}}{\sum }%
\left( \Pi _{j=1}^{n}p_{w_{j}}\right) $

\strut 

$\ \ \ \ \ \ \ \ \ \ \ \ \ \ \ \left( \Pi _{j=1}^{n}\delta
_{(v_{j},x_{j},y_{j}:u_{w_{j}})}\right) \,\,E\left(
L_{w_{1}}^{u_{w_{1}}}...L_{w_{n}}^{u_{w_{n}}}\right) .$

$\square $
\end{proposition}

\strut

From now, rest of this section, we will compute the $D_{G}$-valued cumulants
of the given $D_{G}$-valued random variable $a$. Let $w_{1},...,w_{n}\in
FP(G)$ be finite paths and \ $u_{1},...,u_{n}\in \{1,*\}$. Then, by the M%
\"{o}bius inversion, we have

\strut

(1.13)$\ \ $

\begin{center}
$k_{n}\left( L_{w_{1}}^{u_{1}}~,...,~L_{w_{n}}^{u_{n}}\right) =\underset{\pi
\in NC(n)}{\sum }\widehat{E}(\pi )\left( L_{w_{1}}^{u_{1}}~\otimes
...\otimes ~L_{w_{n}}^{u_{n}}\right) \mu (\pi ,1_{n}),$
\end{center}

\strut

where $\widehat{E}=\left( E^{(n)}\right) _{n=1}^{\infty }$ is the moment
multiplicative bimodule map induced by the conditional expectation $E$ (See
[16]) and where $NC(n)$ is the collection of all noncrossing partition over $%
\{1,...,n\}.$ Notice that if $L_{w_{1}}^{u_{1}}...L_{w_{n}}^{u_{n}}$ does
not have the $*$-axis-property, then

\strut

\begin{center}
$E\left( L_{w_{1}}^{u_{1}}...L_{w_{n}}^{u_{n}}\right) =0_{D_{G}},$
\end{center}

\strut

by Section 2.1. Consider the noncrossing partition $\pi \in NC(n)$ with its
blocks $V_{1},...,V_{k}$. Choose one block $V_{j}=(j_{1},...,j_{k})\in \pi .$
Then we have that

\strut

(1.14) $\ $

\begin{center}
$\widehat{E}(\pi \mid _{V_{j}})\left( L_{w_{1}}^{u_{1}}~\otimes ...\otimes
~L_{w_{n}}^{u_{n}}\right) =E\left(
L_{w_{j_{1}}}^{u_{j_{1}}}d_{j_{2}}L_{w_{j_{2}}}^{u_{j_{2}}}...d_{j_{k}}L_{w_{j_{k}}}^{u_{j_{k}}}\right) , 
$
\end{center}

\strut where

\begin{center}
$d_{j_{i}}=\left\{ 
\begin{array}{ll}
1_{D_{G}} & 
\begin{array}{l}
\text{if there is no inner blocks} \\ 
\text{between }j_{i-1}\text{ \ and }j_{i}\text{ in }V_{j}
\end{array}
\\ 
&  \\ 
L_{v_{j_{i}}}\neq 1_{D_{G}} & 
\begin{array}{l}
\text{if there are inner blocks} \\ 
\text{between }j_{i-1}\text{ \ and }j_{i}\text{ in }V_{j},
\end{array}
\end{array}
\right. $
\end{center}

\strut

where $v_{j_{2}},...,v_{j_{k}}\in V(G).$ So, again by Section 2.1, $\widehat{%
E}(\pi \mid _{V_{j}})\left( L_{w_{1}}^{u_{1}}~\otimes ...\otimes
~L_{w_{n}}^{u_{n}}\right) $ is nonvanishing if and only if $%
L_{w_{j_{1}}}^{u_{j_{1}}}d_{j_{2}}L_{w_{j_{2}}}^{u_{j_{2}}}...d_{j_{k}}L_{w_{j_{k}}}^{u_{j_{k}}} 
$ has the $*$-axis-property, for all \ $j=1,...,n.$

\strut

Assume that

\begin{center}
$\widehat{E}(\pi \mid _{V_{j}})\left( L_{w_{1}}^{u_{1}}~\otimes ...\otimes
~L_{w_{n}}^{u_{n}}\right) =L_{v_{j}}$
\end{center}

and

\begin{center}
$\widehat{E}(\pi \mid _{V_{i}})\left( L_{w_{1}}^{u_{1}}~\otimes ...\otimes
~L_{w_{n}}^{u_{n}}\right) =L_{v_{i}}.$
\end{center}

\strut

If $v_{j}\neq v_{i},$ then the partition-dependent $D_{G}$-moment satisfies
that

\strut

\begin{center}
$\widehat{E}(\pi )\left( L_{w_{1}}^{u_{1}}~\otimes ...\otimes
~L_{w_{n}}^{u_{n}}\right) =0_{D_{G}}.$
\end{center}

\strut

This says that $\widehat{E}(\pi )\left( L_{w_{1}}^{u_{1}}~\otimes ...\otimes
~L_{w_{n}}^{u_{n}}\right) \neq 0_{D_{G}}$ if and only if there exists $v\in
V(G)$ such that

\strut

\begin{center}
$\widehat{E}(\pi \mid _{V_{j}})\left( L_{w_{1}}^{u_{1}}~\otimes ...\otimes
~L_{w_{n}}^{u_{n}}\right) =L_{v},$
\end{center}

\strut

for all \ $j=1,...,k.$

\strut

\begin{definition}
Let $NC(n)$ be the set of all \ noncrossing partition over $\{1,...,n\}$ and
let $L_{w_{1}}^{u_{1}},$ $...,$ $L_{w_{n}}^{u_{n}}\in \left(
W^{*}(G),E\right) $ be $D_{G}$-valued random variables, where $%
u_{1},...,u_{n}\in \{1,*\}.$ We say that the $D_{G}$-valued random variable $%
L_{w_{1}}^{u_{1}}...L_{w_{n}}^{u_{n}}$ is $\pi $-connected if the $\pi $%
-dependent $D_{G}$-moment of it is nonvanishing, for $\pi \in NC(n).$ In
other words, the random variable $L_{w_{1}}^{u_{1}}...L_{w_{n}}^{u_{n}}$ is $%
\pi $-connected, for $\pi \in NC(n),$ if

\strut 

$\ \ \ \ \ \ \ \ \ \ \ \widehat{E}(\pi )\left( L_{w_{1}}^{u_{1}}~\otimes
...\otimes ~L_{w_{n}}^{u_{n}}\right) \neq 0_{D_{G}}.$

\strut 

i.e., there exists a vertex $v\in V(G)$ such that

\strut 

$\ \ \ \ \ \ \ \ \ \ \ \widehat{E}(\pi )\left( L_{w_{1}}^{u_{1}}~\otimes
...\otimes ~L_{w_{n}}^{u_{n}}\right) =L_{v}.$
\end{definition}

\strut

For convenience, we will define the following subset of $NC(n)$ ;

\strut

\begin{definition}
Let $NC(n)$ be the set of all noncrossing partitions over $\{1,...,n\}$ and
fix a $D_{G}$-valued random variable $L_{w_{1}}^{u_{1}}...L_{w_{n}}^{u_{n}}$
in $\left( W^{*}(G),E\right) ,$ where $u_{1},$ ..., $u_{n}\in \{1,*\}.$ For
the fixed $D_{G}$-valued random variable $%
L_{w_{1}}^{u_{1}}...L_{w_{n}}^{u_{n}},$define

\strut 

$\ \ \ C_{w_{1},...,w_{n}}^{u_{1},...,u_{n}}\overset{def}{=}\{\pi \in
NC(n):L_{w_{1}}^{u_{1}}...L_{w_{n}}^{u_{n}}$ is $\pi $-connected$\},$

\strut 

in $NC(n).$ Let $\mu $ be the M\"{o}bius function in the incidence algebra $%
I_{2}.$ Define the number $\mu _{w_{1},...,w_{n}}^{u_{1},...,u_{n}},$ for
the fixed $D_{G}$-valued random variable $%
L_{w_{1}}^{u_{1}}...L_{w_{n}}^{u_{n}},$ by

\strut 

$\ \ \ \ \ \ \ \ \ \ \ \ \ \ \ \mu _{w_{1},...,w_{n}}^{u_{1},...,u_{n}}%
\overset{def}{=}\underset{\pi \in C_{w_{1},...,w_{n}}^{u_{1},...,u_{n}}}{%
\sum }\mu (\pi ,1_{n}).$
\end{definition}

\bigskip \strut

Assume that there exists $\pi \in NC(n)$ such that $%
L_{w_{1}}^{u_{1}}...L_{w_{n}}^{u_{n}}=L_{v}$ is $\pi $-connected. Then $\pi
\in C_{w_{1},...,w_{n}}^{u_{1},...,u_{n}}$ and there exists the maximal
partition $\pi _{0}\in C_{w_{1},...,w_{n}}^{u_{1},...,u_{n}}$ such that $%
L_{w_{1}}^{u_{1}}...L_{w_{n}}^{u_{n}}=L_{v}$ is $\pi _{0}$-connected. Notice
that $1_{n}\in C_{w_{1},...,w_{n}}^{u_{1},...,u_{n}}.$ Therefore, the
maximal partition in $C_{w_{1},...,w_{n}}^{u_{1},...,u_{n}}$ is $1_{n}.$
Hence we have that ;

\bigskip \strut

\begin{lemma}
(Also See [15) Let $L_{w_{1}}^{u_{1}}...L_{w_{n}}^{u_{n}}\in \left(
W^{*}(G),E\right) $ be a $D_{G}$-valued random variable having the $*$%
-axis-property. Then

\strut 

$\ \ \ \ \ \ \ \ \ \ \ E\left( L_{w_{1}}^{u_{1}}...L_{w_{n}}^{u_{n}}~\right)
=\ \widehat{E}(\pi )\left( L_{w_{1}}^{u_{1}}~\otimes ...\otimes
~L_{w_{n}}^{u_{n}}\right) ,$

\strut 

for all $\pi \in C_{w_{1},...,w_{n}}^{u_{1},...,u_{n}}.$ \ $\square $
\end{lemma}

\bigskip \strut \strut \strut

By the previous lemmas, we have that

\strut

\begin{theorem}
(See [15]) Let $n\in 2\mathbb{N}$ and let $%
L_{w_{1}}^{u_{1}},...,L_{w_{n}}^{u_{n}}\in \left( W^{*}(G),E\right) $ be $%
D_{G}$-valued random variables, where $w_{1},...,w_{n}\in FP(G)$ and $%
u_{j}\in \{1,*\},$ $j=1,...,n.$ Then

\strut 

$\ \ \ \ \ \ \ k_{n}\left( L_{w_{1}}^{u_{1}}...L_{w_{n}}^{u_{n}}~\right)
=\mu _{w_{1},...,w_{n}}^{u_{1},...,u_{n}}\cdot
E(L_{w_{1}}^{u_{1}},...,L_{w_{n}}^{u_{n}}),$

\strut 

where $\mu _{w_{1},...,w_{n}}^{u_{1},...,u_{n}}=\underset{\pi \in
C_{w_{1},...,w_{n}}^{u_{1},...,u_{n}}}{\sum }\mu (\pi ,1_{n}).$ \ $\square $
\end{theorem}

\bigskip 

\strut \strut \strut \strut \strut

\strut \strut

\subsection{\strut $D_{G}$-Freeness on $\left( W^{*}(G),E\right) $}

\strut

\strut

Throughout this chapter, let $G$ be a countable directed graph and $\left(
W^{*}(G),E\right) $, the graph $W^{*}$-probability space over its diagonal
subalgebra $D_{G}.$ In this chapter, we will consider the $D_{G}$-valued
freeness of given two random variables in $\left( W^{*}(G),E\right) $. We
will characterize the $D_{G}$-freeness of $D_{G}$-valued random variables $%
L_{w_{1}}$ and $L_{w_{2}},$ where $w_{1}\neq w_{2}\in FP(G).$ And then we
will observe the $D_{G}$-freeness of arbitrary two $D_{G}$-valued random
variables $a_{1}$ and $a_{2}$ in terms of their supports. Let

\strut

(1.15) $\ a=\underset{w\in \mathbb{F}^{+}(G:a),\,u_{w}\in \{1,*\}}{\sum }%
p_{w}L_{w}^{u_{w}}$ \& $b=\underset{w^{\prime }\in \mathbb{F}%
^{+}(G:b),\,u_{w^{\prime }}\in \{1,*\}}{\sum }p_{w^{\prime }}L_{w^{\prime
}}^{u_{w^{\prime }}}$

\strut

be fixed $D_{G}$-valued random variables in $\left( W^{*}(G),E\right) $.

\strut

Now, fix $n\in \mathbb{N}$ and let $\left( a_{i_{1}}^{\varepsilon
_{i_{1}}},...,a_{i_{n}}^{\varepsilon _{i_{n}}}\right) \in
\{a,b,a^{*},b^{*}\}^{n},$ where $\varepsilon _{i_{j}}\in \{1,*\}.$ For
convenience, put

\strut

\begin{center}
$a_{i_{j}}^{\varepsilon _{i_{j}}}=\underset{w_{i_{j}}\in \mathbb{F}%
^{+}(G:a),\,\,u_{j}\in \{1,\ast \}}{\sum }p_{w_{j}}^{(j)}L_{w_{j}}^{u_{j}},$
for \ $j=1,...,n.$
\end{center}

\strut

Then, by the little modification of Section , we have that ;

\strut

(1.9)

\strut

$\ E\left( d_{i_{1}}a_{i_{1}}^{\varepsilon
_{i_{1}}}...d_{i_{n}}a_{i_{n}}^{\varepsilon _{i_{n}}}\right) $

\strut

$\ \ \ \ \ =\,\underset{(v_{i_{1}},...,v_{i_{n}})\in \Pi
_{k=1}^{n}V(G:d_{i_{k}})}{\sum }\left( \Pi _{k=1}^{n}q_{v_{i_{k}}}\right) $

\strut

$\ \ \ \ \ \ \ \ \ \ \ \ \ \ \ \ \ \underset{(w_{i_{1}},...,w_{i_{n}})\in
\Pi _{k=1}^{n}\mathbb{F}^{+}(G:a_{i_{k}}),%
\,w_{i_{j}}=x_{i_{j}}w_{i_{j}}y_{i_{j}},u_{i_{j}}\in \{1,\ast \}}{\sum }%
\left( \Pi _{k=1}^{n}p_{w_{i_{k}}}^{(k)}\right) $

\strut

$\ \ \ \ \ \ \ \ \ \ \ \ \ \ \ \ \ \ \ \ \ \ \ \ \ \ \ \ \ \ \ \ \ \ \left(
\Pi _{j=1}^{n}\delta _{(v_{i_{j}},x_{i_{j}},y_{i_{j}}:u_{i_{j}})}\right)
E\left( L_{w_{i_{1}}}^{u_{i_{1}}}...L_{w_{i_{n}}}^{u_{i_{n}}}\right) .$

\strut

\strut

Therefore, we have that

\strut

(1.10)

\strut

$\ \ k_{n}\left( d_{i_{1}}a_{i_{1}}^{\varepsilon
_{i_{1}}},...,d_{i_{n}}a_{i_{n}}^{\varepsilon _{i_{n}}}\right) =\underset{%
(v_{i_{1}},...,v_{i_{n}})\in \Pi _{k=1}^{n}V(G:d_{i_{k}})}{\sum }\left( \Pi
_{k=1}^{n}q_{v_{i_{k}}}\right) $

\strut

$\ \ \ \ \ \ \ \ \ \ \ \ \ \ \ \ \ \underset{(w_{i_{1}},...,w_{i_{n}})\in
\Pi
_{k=1}^{n}FP(G:a_{i_{k}}),\,w_{i_{j}}=x_{i_{j}}w_{i_{j}}y_{i_{j}},u_{i_{j}}%
\in \{1,\ast \}}{\sum }\left( \Pi _{k=1}^{n}p_{w_{i_{k}}}^{(k)}\right) $

\strut

$\ \ \ \ \ \ \ \ \ \ \ \ \ \ \ \ \ \ \ \ \left( \Pi _{j=1}^{n}\delta
_{(v_{i_{j}},x_{i_{j}},y_{i_{j}}:u_{i_{j}})}\right) \left( \mu
_{w_{i_{1}},...,w_{i_{n}}}^{u_{i_{1}},...,u_{i_{n}}}\cdot E\left(
L_{w_{i_{1}}}^{u_{i_{1}}}...L_{w_{i_{n}}}^{u_{i_{n}}}\right) \right) ,$

\strut

where $\mu _{w_{1},...,w_{n}}^{u_{1},...,u_{n}}=\underset{\pi \in
C_{w_{i_{1}},...,w_{i_{n}}}^{u_{i_{1}},...,u_{i_{n}}}}{\sum }\mu (\pi
,1_{n}) $ and

\bigskip

\begin{center}
$C_{w_{i_{1}},...,w_{i_{n}}}^{u_{i_{1}},...,u_{i_{n}}}=\{\pi \in
NC^{(even)}(n):L_{w_{1}}^{u_{1}}...L_{w_{n}}^{u_{n}}$ is $\pi $-connected$%
\}. $
\end{center}

\strut

So, we have the following proposition, by the straightforward computation ;

\bigskip \strut \strut

\begin{proposition}
Let $a,b\in \left( W^{*}(G),E\right) $ be $D_{G}$-valued random variables,
such that $a\notin W^{*}(\{b\},D_{G}),$ and let $\left(
a_{i_{1}}^{\varepsilon _{i_{1}}},...,a_{i_{n}}^{\varepsilon _{i_{n}}}\right)
\in \{a,b,a^{*},b^{*}\}^{n},$ for $n\in \mathbb{N}\setminus \{1\},$ where $%
\varepsilon _{i_{j}}\in \{1,*\},$ $j=1,...,n.$ Then\strut 

\strut 

\strut (1.11)

\strut 

$\ \ \ \ \ k_{n}\left( d_{i_{1}}a_{i_{1}}^{\varepsilon
_{i_{1}}},...,d_{i_{n}}a_{i_{n}}^{\varepsilon _{i_{n}}}\right) $

\strut 

$\ \ \ \ \ \ =\underset{(v_{1},...,v_{n})=(x,y,...,x,y)\in \Pi
_{j=1}^{n}V(G:d_{j})}{\sum }\left( \Pi _{j=1}^{n}q_{v_{j}}\right) $

\strut 

$\ \ \ \ \ \ \ \ \ \ \underset{(w_{i_{1}},...,w_{i_{n}})\in \left( \Pi
_{k=1}^{n}FP_{*}(G:a_{i_{k}})\right) \cup
W_{*}^{iw_{1},...,w_{n}},\,w_{i_{j}}=x_{i_{j}}w_{i_{j}}y_{i_{j}}}{\sum }%
\left( \Pi _{k=1}^{n}p_{w_{i_{j}}}^{(k)}\right) $

\strut 

$\ \ \ \ \ \ \ \ \ \ \ \ \ \ \left( \Pi _{j=1}^{n}\delta
_{(v_{i_{j}},x_{i_{j}},y_{i_{j}}:u_{i_{j}})}\right) \left( \mu
_{w_{i_{1}},...,w_{i_{n}}}^{u_{i_{1}},...,u_{i_{n}}}\cdot \Pr oj\left(
L_{w_{i_{1}}}^{u_{i_{1}}}...L_{w_{i_{n}}}^{u_{i_{n}}}\right) \right) $

\strut 

\strut where $\mu _{n}=\underset{\pi \in
C_{w_{i_{1}},...,w_{i_{n}}}^{u_{i_{1}},...,u_{i_{n}}}}{\sum }\mu (\pi ,1_{n})
$ and

\strut 

$\ \ \ \ \ \ W_{*}^{w_{1},...,w_{n}}=\{w\in FP_{*}^{c}(G:a)\cup
FP_{*}^{c}(G:b):$

$\ \ \ \ \ \ \ \ \ \ \ \ \ \ \ \ \ \ \ \ \ \ \ \ \ \ \ \ \ \ \ \ \ $both $%
L_{w}^{u_{w}}$ and $L_{w}^{u_{w}\,*}$\ are in $%
L_{w_{1}}^{u_{w_{1}}}...L_{w_{n}}^{u_{w_{n}}}\}.$

$\square $
\end{proposition}

\strut \strut \strut \strut \strut \strut

So, we have the following $D_{G}$-freeness characterization ;

\strut

\begin{corollary}
Let $x$ and $y$ be the $D_{G}$-valued random variables in $\left(
W^{*}(G),E\right) $. The $D_{G}$-valued random variables $a$ and $b$ are
free over $D_{G}$ in $\left( W^{*}(G),E\right) $ if

\strut \strut \strut 

$\ \ \ \ \ \ \ FP_{*}\left( G:P(x,x^{*})\right) \cap FP_{*}\left(
G:Q(y,y^{*})\right) =\emptyset $

and

$\ \ \ \ \ \ \ \ \ \ \ \ \ \ \ W_{*}^{\{P(x,x^{*}),\,Q(y,y^{*})\}}=\emptyset
,$

\strut 

for all $P,Q\in \mathbb{C}[z_{1},z_{2}].$ \ $\square $
\end{corollary}

\strut \strut \strut \strut

We have the above $D_{G}$-freeness characterization, but it is so hard to
use the above characterization. So, we will restrict our interests to the $%
D_{G}$-freeness on the generator set $\{L_{w},L_{w}^{*}:w\in \Bbb{F}%
^{+}(G)\} $ of the graph $W^{*}$-algebra $W^{*}(G).$ In this case, the $%
D_{G} $-freeness on the set is pictorially determined on the given graph $G.$
Now, we will introduce the diagram-distinctness of general finite paths ;

\strut

\begin{definition}
(\textbf{Diagram-Distinctness}) We will say that the finite paths $w_{1}$
and $w_{2}$ are \textbf{diagram-distinct} if $w_{1}$ and $w_{2}$ have
different diagrams in the graph $G.$ Let $X_{1}$ and $X_{2}$ be subsets of $%
FP(G).$ The subsets $X_{1}$ and $X_{2}$ are said to be diagram-distinct if $%
x_{1}$ and $x_{2}$ are diagram-distinct for all pairs $(x_{1},x_{2})$ $\in $ 
$X_{1}\times X_{2}.$
\end{definition}

\strut \strut

Let $H$ be a directed graph with $V(H)=\{v_{1},v_{2}\}$ and $%
E(H)=\{e_{1}=v_{1}e_{1}v_{2},e_{2}=v_{2}e_{2}v_{1}\}.$ Then $l=e_{1}e_{2}$
is a loop in $FP(H)$ (i.e., $l\in loop(H)$). Moreover, it is a basic loop
(i.e., $l\in Loop(H)$). However, if we have a loop $%
w=e_{1}e_{2}e_{1}e_{2}=l^{2},$ then it is not a basic loop. i.e.,

\strut

\begin{center}
$l^{2}\in loop(H)\,\setminus \,Loop(H).$
\end{center}

\strut

If the graph $G$ contains at least one basic loop $l\in FP(G)$, then we have

\strut

\begin{center}
$\{l^{n}:n\in \mathbb{N}\}\subseteq loop(G)$ \ and \ $\{l\}\subseteq
Loop(G). $
\end{center}

\strut

Suppose that $l_{1}$ and $l_{2}$ are not diagram-distinct. Then, by
definition, there exists $w\in Loop(G)$ such that $l_{1}=w^{k_{1}}$ and $%
l_{2}=w^{k_{2}},$ for some $k_{1},k_{2}\in \mathbb{N}.$ On the graph $G,$
indeed, $l_{1}$ and $l_{2}$ make the same diagram. On the other hands, we
can see that if $w_{1}\neq w_{2}\in loop^{c}(G),$ then they are
automatically diagram-distinct. In [15], we found the $D_{G}$-freeness
characterization on the generator set of $W^{*}(G),$ as follows ;

\strut \strut \strut

\begin{theorem}
(See [15]) Let $w_{1},w_{2}\in FP(G)$ be finite paths. The $D_{G}$-valued
random variables $L_{w_{1}}$ and $L_{w_{2}}$ in $\left( W^{*}(G),E\right) $
are free over $D_{G}$ if and only if $w_{1}$ and $w_{2}$ are
diagram-distinct. $\square $
\end{theorem}

\strut \strut \strut

Let $a$ and $b$ be the given $D_{G}$-valued random variables. We can get the
necessary condition for the $D_{G}$-freeness of $a$ and $b,$ in terms of
their supports. \strut \strut \strut Recall that we say that the two subsets 
$X_{1}$ and $X_{2}$ of $FP(G)$ are said to be diagram-distinct if $x_{1}$
and $x_{2}$ are diagram-distinct, for all pairs $(x_{1},x_{2})$ $\in $ $X_{1}
$ $\times $ $X_{2}.$

\strut \strut

\begin{proposition}
(See [15]) Let $a,b\in \left( W^{*}(G),E\right) $ be $D_{G}$-valued random
variables with their supports $\mathbb{F}^{+}(G:a)$ and $\mathbb{F}^{+}(G:b).
$ The $D_{G}$-valued random variables $a$ and $b$ are free over $D_{G}$ in $%
\left( W^{*}(G),E\right) $ if $FP(G:a_{1})$ and $FP(G:a_{2})$ are
diagram-distinct. $\square $
\end{proposition}

\strut

\strut \strut

\strut

\subsection{$D_{G}$\protect\bigskip -valued Random Variables}

\bigskip

\strut In this section, we will consider certain $D_{G}$-valued random
variables. Let $G$ be a countable directed graph and let $\left(
W^{*}(G),E\right) $ be the graph $W^{*}$-probability space over its diagonal
subalgebra $D_{G}.$ In [17], we showed that the $D_{G}$-semicircularity, the 
$D_{G}$-evenness and the $D_{G}$-valued R-diagonality of $D_{G}$-valued
random variables can be characterized by the graphical expression on the
graph $G.$

\strut

Let $B$ be a von Neumann algebra and $A,$ a von Neumann algebra over $B$ and
let $F:A\rightarrow B$ be a conditional expectation. Then we have a $W^{*}$%
-probability space $(A,F)$ over $B.$ The $B$-valued random variable $a\in
(A,F)$ is called a $B$-semicircular element if it is self-adjoint and the
only nonvanishing $B$-cumulant of $a$ is the second one. i.e., a $B$%
-semicircular element $a$ satisfies that

\strut

\begin{center}
$k_{n}^{F}(a,...,a)=\left\{ 
\begin{array}{ll}
k_{2}^{F}(a,a)\neq 0_{B} & \text{if }n=2 \\ 
&  \\ 
0_{B} & \text{otherwise,}
\end{array}
\right. $
\end{center}

\strut

where $k_{n}^{F}(..)$ is the $B$-cumulant bimodule map induced by the
conditional expectation $F.$ Suppose $x$ is a $B$-valued random variable in $%
\left( A,F\right) .$ We say that the random variable $x$ is $B$-even if it
is self-adjoint and all odd $B$-cumulants vanish. Equivalently, the
self-adjoint operator $x$ is $B$-even if all odd $B$-moments vanish. Now,
assume that $y\in \left( A,F\right) $ is a $B$-valued random variable. If
the only nonvanishing mixed $B$-cumulant of $y$ and $y^{*}$ are alternating $%
B$-cumulants, i.e., if the only nonvanishing $B$-cumulants of $y$ and $y^{*}$
are

\strut

\begin{center}
$k_{2n}\left( y,y^{*},...,y,y^{*}\right) $ \ and \ $k_{2n}\left(
y^{*},y,...,y^{*},y\right) ,$
\end{center}

\strut

for all $n\in \Bbb{N},$ then the $D_{G}$-valued random variable $y$ (and $%
y^{*}$) is called the $B$-valued R-diagonal.

\strut

The $B$-semicircular elements, $B$-even elements and $B$-valued R-diagonal
elements play important role in Free Probability. The following theorem
shows that the graph $W^{*}$-probability spaces contain such random
variables. So, the graph $W^{*}$-probability spaces contain rich free
probabilistic objects.

\strut \strut

\begin{proposition}
(See [17]) Let $w\in \Bbb{F}^{+}(G)$ and let $L_{w}$ be the corresponding $%
D_{G}$-valued random variable in $\left( W^{*}(G),E\right) .$ Then

\strut

(1) if $w$ is a loop, then $L_{w}+L_{w}^{*}$ is $D_{G}$-semicircular.

\strut

(2) if $w$ is a finite path, then $L_{w}+L_{w}^{*}$ is $D_{G}$-even.

\strut

(3) if $w$ is a finite path, then $L_{w}$ and $L_{w}^{*}$ are $D_{G}$-valued
R-diagonal. $\square $
\end{proposition}

\strut \strut

\strut \strut

\strut \strut

\section{Vertex Compressed Graph $W^{*}$-Probability}

\strut

\strut

\strut

In this chapter, we will consider the vertex-compressed graph $W^{*}$%
-probability. Throughout this chapter, let $G$ be a countable directed graph
and let $\left( W^{*}(G),E\right) $ be the graph $W^{*}$-probability space
over the diagonal subalgebra $D_{G}.$ Let $v_{0}\in V(G)$ be a vertex and
let's fix this vertex. Then we can take the projection $L_{v_{0}}\in \left(
W^{*}(G),E\right) .$ We will consider the compressed $W^{*}$-algebra $%
L_{v_{0}}(W^{*}(G))L_{v_{0}}$ by $L_{v_{0}}$ and observe the compressed
probability space on $\left( L_{v_{0}}W^{*}(G)L_{v_{0}},E_{v_{0}}\right) ,$
where

\strut

\begin{center}
$E_{v_{0}}:L_{v_{0}}W^{*}(G)L_{v_{0}}\rightarrow D_{G}$
\end{center}

\strut

is the compressed conditional expectation defined by

$\strut $

\begin{center}
$E_{v_{0}}\overset{def}{=}E\mid _{L_{v_{0}}W^{*}(G)L_{v_{0}}}$.
\end{center}

\strut

It is easy to see that the $v_{0}$-compressed conditional expectation $%
E_{v_{0}}$ can be regarded as a linear functional from the $v_{0}$%
-compressed graph $W^{*}$-algebra $L_{v_{0}}W^{*}(G)L_{v_{0}}$ onto $%
\mathbb{C}=\mathbb{C}\xi _{v_{0}}=L_{v_{0}}D_{G}L_{v_{0}}.$ Indeed, let

\strut

\begin{center}
$a=\underset{w\in \Bbb{F}^{+}(G:a),\,u_{w}\in \{1,*\}}{\sum }%
p_{w}L_{w}^{u_{w}}\in \left( W^{*}(G),E\right) .$
\end{center}

\strut

Then, for any summand $L_{w}^{u_{w}}$ of $a,$ $w\in \Bbb{F}^{+}(G:a),$

\strut

\begin{center}
$L_{v_{0}}L_{w}^{u_{w}}L_{v_{0}}=\left\{ 
\begin{array}{lll}
L_{v_{0}} &  & \text{if }w=v_{0} \\ 
L_{w}^{u_{w}} &  & \text{if }l=v_{0}lv_{0} \\ 
0_{D_{G}} &  & \text{otherwise,}
\end{array}
\right. $
\end{center}

\strut

for all $u_{w}=1,*.$ Thus $E_{v_{0}}$ maps $L_{v_{0}}W^{*}(G)L_{v_{0}}$
linearly onto $\Bbb{C}\xi _{v_{0}}.$\strut This shows that when we want to
compute the $L_{v_{0}}D_{G}L_{v_{0}}$-valued moments and cumulants of a
compressed random variable $L_{v_{0}}aL_{v_{0}}$, the \textbf{trivial}
moments and cumulants contain the full free probabilistic information of the
compressed random variable $L_{v_{0}}aL_{v_{0}}.$ Also, we can easily verify
that

\strut

\begin{center}
$E_{v_{0}}(x)=\,<\xi _{v_{0}},x\xi _{v_{0}}>\,\in \mathbb{C\cdot \xi }%
_{v_{0}},$ \ $\forall x\in L_{v_{0}}W^{\ast }(G)L_{v_{0}}.$
\end{center}

\strut

Again, remark that

\begin{center}
$L_{v_{0}}aL_{v_{0}}=p_{v_{0}}[L_{v_{0}}]+\underset{w=v_{0}wv_{0}\in
loop(G:a),\,u_{w}\in \{1,*\}}{\sum }p_{w}L_{w}^{u_{w}},$
\end{center}

\strut

where $[L_{v_{0}}]=L_{v_{0}}$ if $v_{0}\in V(G:a)$ and $%
[L_{v_{0}}]=0_{D_{G}},$ otherwise. Therefore, we can get that

\strut

\begin{center}
$E_{v_{0}}\left( L_{v_{0}}aL_{v_{0}}\right) =p_{v_{0}}[L_{v_{0}}]$ $\in
D_{G}.$
\end{center}

\strut

Hence, we can consider the $v_{0}$-compressed graph $W^{*}$-probability
space over $D_{G},$ $\left( L_{v_{0}}W^{*}(G)L_{v_{0}},E_{v_{0}}\right) $,
as a (scalar-valued) $W^{*}$-probability space. Let's regard $%
L_{v_{0}}D_{G}L_{v_{0}}$ as $\mathbb{C}.$

\strut

\strut

\subsection{Vertex Compressed Graph $W^{*}$-Probability Spaces}

\strut

\strut

In this section, we will consider the vertex compressed graph $W^{*}$%
-probability space as a (scalar-valued) $W^{*}$-probability space (over $%
\mathbb{C}$). Let $a\in \left( W^{*}(G),E\right) $ be a $D_{G}$-valued
random variable with $FP(G:a)$ $\subseteq $ $loop(G).$ Then such $D_{G}$%
-valued random variable $a$ is called a loop operator in $\left(
W^{*}(G),E\right) .$ By definition, the $v_{0}$-compressed random variable $%
x $ has its form of

\strut

\begin{center}
$x=p_{v_{0}}[L_{v_{0}}]+\underset{w\in loop_{v_{0}}(G:a),\,u_{w}\in \{1,*\}}{%
\sum }p_{w}L_{w}^{u_{w}}.$
\end{center}

\strut

So, every $v_{0}$-compressed random variables is a loop operator, where

\strut

\begin{center}
$loop_{v_{0}}(G:a)\overset{def}{=}\{l\in loop(G:a):l=v_{0}lv_{0}\}.$
\end{center}

\strut

\begin{definition}
Let $G$ be a countable directed graph and fix $v_{0}\in V(G).$ Define the
compressed $W^{*}$-algebra, $W_{v_{0}}^{*}(G)$ $\overset{denote}{=}$ $%
L_{v_{0}}W^{*}(G)L_{v_{0}}$ and we will call it the $v_{0}$-compressed graph 
$W^{*}$-algebra. Now, define the linear functional $E_{v_{0}}:$ $%
W_{v_{0}}^{*}(G)\rightarrow \mathbb{C}$ by

\strut

\ \ \ $\ \ \ \ \ \ E_{v_{0}}\left( x\right) =\,<\xi _{v_{0}},\,x\xi
_{v_{0}}>,$ \ for all \ $x\in W^{*}(G).$

\strut

We will call the algebraic pair $\left( W_{v_{0}}^{*}(G),E_{v_{0}}\right) ,$
the $v_{0}$-compressed graph $W^{*}$-probability space. Let $D_{G}$ be the
diagonal subalgebra. Denote $L_{v_{0}}D_{G}L_{v_{0}}$ by $D_{G}^{v_{0}},$
for convenience.
\end{definition}

\strut \strut

Notice that we can regard the $v_{0}$-compressed graph $W^{*}$-probability
space $\left( W_{v_{0}}^{*}(G),E_{v_{0}}\right) $ as a scalar-valued $W^{*}$%
-probability space. So, we can define the $n$-th moments and $n$-th
cumulants, in the sense of Nica and Speicher (See [1] and [19]). In our
notation, they are just \textbf{trivial} $D_{G}^{v_{0}}$-valued moments and
cumulants ;

\strut

\begin{definition}
Let $\left( W_{v_{0}}^{*}(G),E_{v_{0}}\right) $ be a $v_{0}$-compressed
graph $W^{*}$-probability space and let $a_{v_{0}}=L_{v_{0}}aL_{v_{0}}\in
\left( W_{v_{0}}^{*}(G),E_{v_{0}}\right) $ be a random variable. Then the $n$%
-th moment of $a_{v_{0}}$ is $E_{v_{0}}\left( a_{v_{0}}^{n}\right) $ and the 
$n$-th cumulant of $a_{v_{0}}$ is $k_{n}^{E_{v_{0}}}\left(
a_{v_{0}},...,a_{v_{0}}\right) .$ If there is no confusion, we will denote
the $n$-th compressed cumulant $k_{n}^{E_{v_{0}}}\left(
a_{v_{0}},...,a_{v_{0}}\right) ,$ by $k_{n}\left(
a_{v_{0}},...,a_{v_{0}}\right) ,$ for all $n\in \mathbb{N}.$ Define the $%
v_{0}$-compressed moment series of $a_{v_{0}}$ by

\strut

$\ \ \ \ \ \ \ \ \ \ \ \ \ M_{a_{v_{0}}}^{(v_{0})}(z)=\sum_{n=1}^{\infty
}E_{v_{0}}\left( a_{v_{0}}^{n}\right) \,z^{n}\in \Bbb{C}[z],$

\strut

and, define the $v_{0}$-compressed R-transform of $a_{v_{0}}$ by

\strut

$\ \ \ \ \ \ \ \ \ \ \ \ \ R_{a_{v_{0}}}^{(v_{0})}(z)=\sum_{n=1}^{\infty
}k_{n}\left( a_{v_{0}},...,a_{v_{0}}\right) \,z^{n}\in \Bbb{C}[z]$.
\end{definition}

\strut

Observe that if $a\in W^{*}(G)$ is given as above, then, by Section 1.3,

\strut

$\ \ \
(L_{v_{0}}aL_{v_{0}})^{n}=(L_{v_{0}}aL_{v_{0}})(L_{v_{0}}aL_{v_{0}})...(L_{v_{0}}aL_{v_{0}}) 
$

\strut

$\ \ \ \ \ \ \ \ \ \
=L_{v_{0}}aL_{v_{0}}a...L_{v_{0}}aL_{v_{0}}=L_{v_{0}}aL_{v_{0}}^{2}aL_{v_{0}}^{2}a...L_{v_{0}}^{2}aL_{v_{0}} 
$

\bigskip

$\ \ \ \ \ \ \ \ \ =\left( L_{v_{0}}aL_{v_{0}}\right) \left(
L_{v_{0}}aL_{v_{0}}\right) ...\left( L_{v_{0}}aL_{v_{0}}\right) $

\strut

(2.1.1)

\strut

$\ \ \ \ \ \ \ \ \ \ =\underset{(w_{1},...,w_{n})\in \mathbb{F}%
^{+}(G:a)^{n},\,w_{j}=v_{0}w_{j}v_{0},\,u_{w_{j}}\in \{1,\ast \}}{\sum }%
\left( \Pi _{j=1}^{n}p_{w_{j}}\right)
L_{w_{1}}^{u_{w_{1}}}...L_{w_{n}}^{u_{w_{n}}}.$

\strut

One may be tempted to use $loop_{v_{0}}(G:a),$ instead of using $\mathbb{F}%
^{+}(G:a),$ in the formula (2.1.1). However, we need to consider the case
when $v_{0}\in V(G:a),$ in general. (Clearly, if $v_{0}\in V(G:a),$ then $%
v_{0}=v_{0}v_{0}v_{0}.$) That's why we used $\mathbb{F}^{+}(G:a),$ in
(2.1.1).

\strut \strut \strut

\strut

\strut

\subsection{Vertex Compressed Moments and Cumulants}

\strut

\strut

Throughout this section, let $G$ be a countable directed graph and let $%
\left( W^{*}(G),E\right) $ be the graph $W^{*}$-probability space over the
diagonal subalgebra $D_{G}.$ Let $v_{0}\in V(G)$ be the fixed vertex. Then
we can construct the $v_{0}$-compressed graph $W^{*}$-probability space $%
\left( W_{v_{0}}^{*}(G),E_{v_{0}}\right) ,$ as a scalar-valued (or $%
D_{G}^{v_{0}}$-valued) $W^{*}$-probability space with the linear functional $%
E_{v_{0}}:W_{v_{0}}^{*}(G)\rightarrow D_{G}^{v_{0}}\simeq \mathbb{C}.$

\strut

Let's consider the $n$-th cumulants of certain $v_{0}$-compressed random
variables in $\left( W_{v_{0}}^{*}(G),\,E_{v_{0}}\right) $ ;\strut

\strut \strut

\begin{proposition}
Let $\left( W_{v_{0}}^{*}(G),E_{v_{0}}\right) $ be the $v_{0}$-compressed
graph $W^{*}$-probability space. Let $a_{v_{0}}=L_{v_{0}}aL_{v_{0}}\in
\left( W_{v_{0}}^{*}(G),E_{v_{0}}\right) $ be a random variable, where $a\in
\left( W^{*}(G),E\right) $ and assume that $loop(G:a)$ consists of all
mutually \textbf{diagram-distinct} loops. Then the $n$-th (scalar-valued)
cumulants of $a_{v_{0}}$ are

\strut

$\ \ k_{n}^{(E_{v_{0}})}\left( \underset{n-times}{\underbrace{%
a_{v_{0}},.....,a_{v_{0}}}}\right) =\left\{ 
\begin{array}{ll}
\underset{l\in loop_{*}^{v_{0}}(G:a)}{\sum }2\left( p_{l}p_{l^{t}}\right)
^{2}L_{v_{0}} & \text{if }n=2 \\ 
&  \\ 
0 & \text{otherwise,}
\end{array}
\right. $

\strut

where $loop_{*}^{v_{0}}(G:a)=\{l\in FP_{*}(G:a):l=v_{0}lv_{0}\}.$
\end{proposition}

\strut

\begin{proof}
Suppose that $a_{v_{0}}=L_{v_{0}}aL_{v_{0}}$ is the $v_{0}$-compressed
random variable in $(W_{v_{0}}^{\ast }(G),E_{v_{0}})\subset (W^{\ast
}(G),E). $ As we have seen before, we have that

\strut

$\ a_{v_{0}}=[p_{v_{0}}L_{v_{0}}]+\underset{l\in loop_{*}^{v_{0}}(G:a)}{\sum 
}\left( p_{l}L_{l}+p_{l^{t}}L_{l}^{*}\right) +\underset{w=v_{0}wv_{0}\in
loop_{*}^{c}(G:a),~u_{w}\in \{1,*\}}{\sum }p_{w}L_{w}^{u_{w}},$

\strut

denoted by $[p_{v_{0}}L_{v_{0}}]+a_{(\ast )}^{v_{0}}+a_{(non-\ast
)}^{v_{0}}. $ Then

\strut

$\ k_{n}^{E_{v_{0}}}\left( a_{v_{0}},...,a_{v_{0}}\right) =k_{n}\left(
a_{(\ast )}^{v_{0}},...,a_{(\ast )}^{v_{0}}\right) $

\strut

$\ \ \ =k_{n}\left( \underset{l\in loop_{\ast }^{v_{0}}(G:a)}{\sum }\left(
p_{l}L_{l}+p_{l^{t}}L_{l}^{\ast }\right) ,...,\underset{l\in loop_{\ast
}^{v_{0}}(G:a)}{\sum }\left( p_{l}L_{l}+p_{l^{t}}L_{l}^{\ast }\right)
\right) $

\strut

$\ \ \ =\underset{l\in loop_{\ast }^{v_{0}}(G:a)}{\sum }~k_{n}\left(
(p_{l}L_{l}+p_{l^{t}}L_{l}^{\ast })~,...,~(p_{l}L_{l}+p_{l^{t}}L_{l}^{\ast
})\right) $

\strut

since $p_{l_{1}}L_{l_{1}}+p_{l_{1}^{t}}L_{l_{1}}^{*}$ and $%
p_{l_{2}}L_{l_{2}}+p_{l_{2}^{t}}L_{l_{2}}^{*}$ are free in $%
(W_{v_{0}}^{*}(G),E_{v_{0}}),$ whenever $l_{1}\neq l_{2}$ in $%
loop_{*}^{v_{0}}(G:a),$ by assumption.

\strut

$\ \ \ =\underset{l\in loop_{*}^{v_{0}}(G:a)}{\sum }\,\underset{%
(u_{1},...,u_{n})\in \{1,*\}^{n}}{\sum }k_{n}\left(
p_{l_{u_{1}}}L_{l}^{u_{1}},...,p_{l_{u_{n}}}L_{l}^{u_{n}}\right) $

\strut

$\ \ =\underset{l\in loop_{*}^{v_{0}}(G:a)}{\sum }\,\underset{%
(u_{1},...,u_{n})\in \{1,*\}^{n}}{\sum }\left( \Pi
_{j=1}^{n}p_{l_{u_{j}}}\right) $

\strut

$\ \ \ \ \ \ \ \ \ \ \ \ \ \ \ \ \ \ \ \ \ \ \ \ \ k_{n}\left(
L_{l}^{u_{1}},...,L_{l}^{u_{n}}\right) $

\strut

by the bilinearity of $k_{n}(...),$ where $p_{l_{u_{j}}}=p_{l}$ if $u_{j}=1$
and $p_{l_{u_{j}}}=p_{l^{t}}$ if $u_{j}=*$

\strut

$\ \ =\left\{ 
\begin{array}{ll}
\underset{l\in loop_{*}^{v_{0}}(G:a)}{\sum }\left( p_{l}p_{l^{t}}\right)
^{2}k_{2}\left( L_{l}+L_{l}^{*},L_{l}+L_{l}^{*}\right) & \text{if }n=2 \\ 
&  \\ 
0 & \text{otherwise}
\end{array}
\right. $

\strut

$\ \ =\left\{ 
\begin{array}{ll}
\underset{l\in loop_{*}^{v_{0}}(G:a)}{\sum }2\left( p_{l}p_{l^{t}}\right)
^{2}L_{v_{0}} & \text{if }n=2 \\ 
&  \\ 
0 & \text{otherwise,}
\end{array}
\right. $

\strut

by the $D_{G}$-semicircularity of $L_{l}+L_{l}^{*}.$
\end{proof}

\strut \strut \strut

In the following theorem, we will compute the moments of an arbitrary $v_{0}$%
-compressed random variable ;

\strut

\begin{theorem}
Let $\left( W_{v_{0}}^{*}(G),E_{v_{0}}\right) $ be the $v_{0}$-compressed
graph $W^{*}$-probability space, as a scalar-valued $W^{*}$-probability
space. Let $a\in W^{*}(G)$ and let $a_{v_{0}}=L_{v_{0}}aL_{v_{0}}\in
W_{v_{0}}^{*}(G).$ Then

\strut

\ \ \ $\ \ \ E_{v_{0}}\left( a_{v_{0}}^{n}\right) =\underset{%
(w_{1},...,w_{n})\in loop_{*}^{v_{0}}(G:a),~u_{j}\in \{1,*\}}{\sum }\left(
\Pi _{j=1}^{n}p_{w_{j}}\right) \cdot L_{v_{0}}$

\strut

for all $n\in \mathbb{N},$ in $\mathbb{C}\xi _{v_{0}},$ where

\strut

\ \ \ \ \ $\ \ \ loop_{*}^{v_{0}}(G:a)=\{l\in FP_{*}(G:a):l=v_{0}lv_{0}\}.$

\strut

In particular, in this case, $loop_{*}^{v_{0}}(G:a)=FP_{*}(G:a_{v_{0}}).$
\end{theorem}

\strut

\begin{proof}
Consider the $v_{0}$-compressed random variable $%
a_{v_{0}}=L_{v_{0}}aL_{v_{0}},$ for the fixed $D_{G}$-valued random
variable, $a=a_{d}+a_{(\ast )}+a_{(non-\ast )}\in (W^{\ast }(G),E).$ Then we
have that

\strut

$\ \ a_{v_{0}}=L_{v_{0}}aL_{v_{0}}$

$\strut $

$\ \ \ \ \ \ \ =L_{v_{0}}a_{d}L_{v_{0}}+L_{v_{0}}a_{(\ast
)}L_{v_{0}}+L_{v_{0}}a_{(non-\ast )}L_{v_{0}}$

\strut

$\ \ \ \ \ \ \ =[p_{v_{0}}L_{v_{0}}]+\underset{l=v_{0}lv_{0}\in loop_{\ast
}^{v_{0}}(G:a)}{\sum }\left( p_{l}L_{l}+p_{l^{t}}L_{l}^{\ast }\right) $

$\ \ \ \ \ \ \ \ \ \ \ \ \ \ \ \ \ \ \ \ \ \ \ \ \ \ \ \ \ \ \ \ \ \ \ +%
\underset{w=v_{0}wv_{0}\in loop_{\ast }^{c}(G:a),~u_{w}\in \{1,\ast \}}{\sum 
}p_{w}L_{w}^{u_{w}},$

\strut

where $[p_{v_{0}}L_{v_{0}}]=p_{v_{0}}L_{v_{0}}$ if $v_{0}\in V(G:a)$ and $%
[p_{v_{0}}L_{v_{0}}]=0$ if $v_{0}\notin V(G:a).$ So, this $v_{0}$-compressed
random variable $a_{v_{0}}$ is an addition of $[p_{v_{0}}L_{v_{0}}]$ and the
loop operator $L_{v_{0}}a_{(*)}L_{v_{0}}+L_{v_{0}}a_{(non-*)}L_{v_{0}},$ in $%
W^{*}(G),$ centered at $v_{0}\in V(G).$ Therefore, by Section 1.3, we have
that

\strut

$\ \ \ \ \ \ \ \ \ E(a_{v_{0}}^{n})=\underset{(w_{1},...,w_{n})\in
loop_{*}^{v_{0}}(G:a),~u_{j}\in \{1,*\}}{\sum }\left( \Pi
_{j=1}^{n}p_{w_{j}}\right) \cdot L_{v_{0}}.$

\strut \strut
\end{proof}

\strut \strut \strut \bigskip \strut

Recall that in Section 1.3, the $D_{G}$-valued cumulants of the fixed random
variable $a$ is easily gotten by multiplying

$\strut $

\begin{center}
$\mu _{w_{1},...,w_{n}}^{u_{w_{1}},...,u_{w_{n}}}=\underset{\pi \in
C_{w_{1},...,w_{n}}^{u_{w_{1}},...,u_{w_{n}}}}{\sum }\mu (\pi ,1_{n})$
\end{center}

\strut

to each summand $E\left(
L_{w_{1}}^{u_{w_{1}}}...L_{w_{n}}^{u_{w_{n}}}\right) $ of the $D_{G}$-valued
moments of $a.$ This happens because of the $*$-axis-property. So, we can
get the following cumulants of arbitrary $v_{0}$-compressed random
variables, by the previous theorem ;

\strut \strut

\begin{corollary}
Let $\left( W_{v_{0}}^{*}(G),E_{v_{0}}\right) $ be the $v_{0}$-compressed
graph $W^{*}$-probability space, as a scalar-valued $W^{*}$-probability
space. Let $a\in W^{*}(G)$ and let $a_{v_{0}}=L_{v_{0}}aL_{v_{0}}\in
W_{v_{0}}^{*}(G).$ Then

\strut

$k_{n}^{(E_{v_{0}})}=\underset{(w_{1},...,w_{n})\in
loop_{*}^{v_{0}}(G:a),~u_{j}\in \{1,*\}}{\sum }\left( \Pi
_{j=1}^{n}p_{w_{j}}\right) \cdot \mu
_{w_{1},...,w_{n}}^{u_{w_{1}},...,u_{w_{n}}}L_{v_{0}},$

\strut

for all $n\in \Bbb{N}.$ \ $\square $
\end{corollary}

\strut

\strut

\strut \strut

\subsection{Vertex-Compressed Freeness}

\strut

\strut

In this section, we will consider the vertex-compressed freeness on the
graph $W^{*}$-probability space over the diagonal subalgebra $D_{G}.$ Let $%
\left( W^{*}(G),E\right) $ be a graph $W^{*}$-probability space over the
diagonal subalgebra $D_{G}$ and let $\left(
W_{v_{0}}^{*}(G),E_{v_{0}}\right) $ be the $v_{0}$-compressed graph $W^{*}$%
-probability space as a scalar-valued $W^{*}$-probability space (over $%
D_{G}^{v_{0}}=L_{v_{0}}D_{G}L_{v_{0}}\simeq \mathbb{C}$).

\strut

Let $X$ and $Y$ be subalgebras of the $v_{0}$-compressed $W^{*}$-algebra, $%
W_{v_{0}}^{*}(G).$ We say that the subalgebras $X$ and $Y$ are free if all
mixed cumulants of $X$ and $Y$ vanish (with respect to the $v_{0}$%
-compressed conditional expectation or linear functional $%
E_{v_{0}}:W_{v_{0}}^{*}(G)\rightarrow D_{G}^{v_{0}}\simeq \mathbb{C}$).
Also, two random variables $x$ and $y$ are free in $\left(
W_{v_{0}}^{*}(G),E_{v_{0}}\right) $ if $x\notin W^{*}(\{y\}),$ $y\notin
W^{*}(\{x\})$ and if all mixed cumulants of $x$ and $y$ vanish. Suppose that 
$x$ and $y$ are random variables in $\left(
W_{v_{0}}^{*}(G),E_{v_{0}}\right) .$ Then there exists operators $a$ and $b$
in $W^{*}(G)$ such that

$\strut $

\begin{center}
$x=L_{v_{0}}aL_{v_{0}}$ \ \ \ \ and \ \ \ \ $y=L_{v_{0}}bL_{v_{0}}.$
\end{center}

\strut

Recall that the $D_{G}$-valued random variables

\strut

\begin{center}
$a_{1}=\underset{v_{1}\in V(G:a_{1})}{\sum }p_{v_{1}}^{(1)}L_{v_{1}}+%
\underset{w_{1}\in FP(G:a_{1}),\,u_{w_{1}}\in \{1,*\}}{\sum }%
p_{w_{1}}^{(1)}L_{w_{1}}^{u_{w_{1}}}$
\end{center}

and

\begin{center}
$a_{2}=\underset{v_{2}\in V(G:a_{2})}{\sum }p_{v_{2}}^{(2)}L_{v_{2}}+%
\underset{w_{2}\in FP(G:a_{2}),\,u_{w_{2}}\in \{1,*\}}{\sum }%
p_{w_{2}}^{(2)}L_{w_{2}}^{u_{w_{2}}}$
\end{center}

\strut

are free over $D_{G}$ if $FP(G:a_{1})$ and $FP(G:a_{2})$ are
diagram-distinct.

\strut

By the above result, we have that ;

\bigskip \strut \strut

\begin{theorem}
Let $a_{1}$ and $a_{2}$ be $D_{G}$-valued random variables in $\left(
W^{*}(G),E\right) $ and assume $FP(G:a_{1})$ and $FP(G:a_{2})$ are
diagram-distinct$.$ Then the $v_{0}$-compressed random variables $%
x=L_{v_{0}}a_{1}L_{v_{0}}$ and $y=L_{v_{0}}a_{2}L_{v_{0}}$ are free in $%
\left( W_{v_{0}}^{*}(G),E_{v_{0}}\right) .$
\end{theorem}

\strut

\begin{proof}
We have that

\strut

$\ \ \ \ \ \ \ x=p_{v_{0}}^{(1)}[L_{v_{0}}]+\underset{l_{1}\in
loop^{v_{0}}(G:a_{1}),\,u_{l_{1}}\in \{1,*\}}{\sum }%
p_{l_{1}}^{(1)}L_{l_{1}}^{u_{l_{1}}}$

and

$\ \ \ \ \ \ \ y=p_{v_{0}}^{(2)}[L_{v_{0}}]+\underset{l_{2}\in
loop^{v_{0}}(G:a_{2}),\,u_{l_{2}}\in \{1,*\}}{\sum }%
p_{l_{2}}^{(2)}L_{l_{2}}^{u_{l_{2}}}.$

\strut

Since $FP(G:a_{1})$ and $FP(G:a_{2})$ are diagram-distinct, $%
loop_{v_{0}}(G:a_{1})$ and $loop_{v_{0}}(G:a_{2})$ are diagram-distinct.
Notice that

\strut

$\ \ \ \ \ \ \ \ \ \ \ \ \ \ \ FP(G:x)=loop_{v_{0}}(G:a_{1})$

and

$\ \ \ \ \ \ \ \ \ \ \ \ \ \ \ FP(G:y)=loop_{v_{0}}(G:a_{2}).$

\strut

Thus, the $v_{0}$-compressed random variables $x$ and $y$ are free over $%
D_{G}$ in $\left( W^{*}(G),E\right) .$ Remark that the $v_{0}$-compressed
random variables $x$ and $y$ are scalar-valued random variables in $\left(
W_{v_{0}}^{*}(G),E_{v_{0}}\right) $ and the compressed moments and cumulants
are same as the $D_{G}$-valued moments and cumulants of $x$ and $y$ over $%
D_{G}^{v_{0}}\hookrightarrow D_{G}.$ Therefore, the $v_{0}$-compressed
random variables $x$ and $y$ are free in $\left(
W_{v_{0}}^{*}(G),E_{v_{0}}\right) .$
\end{proof}

\strut \strut \strut

We also have the following general case. This shows that the compressed
freeness is preserved by the $D_{G}$-freeness.

\strut

\begin{theorem}
Let $a$ and $b$ be $D_{G}$-valued random variables in the graph $W^{*}$%
-probability space $\left( W^{*}(G),E\right) $ over the diagonal subalgebra $%
D_{G}.$ If they are free over $D_{G}$ in $\left( W^{*}(G),E\right) ,$ then
the corresponding $v_{0}$-compressed random variables $x=L_{v_{0}}aL_{v_{0}}$
and $\ y=L_{v_{0}}bL_{v_{0}}$ are free in $\left(
W_{v_{0}}^{*}(G),E_{v_{0}}\right) .$
\end{theorem}

\strut

\begin{proof}
Suppose that there exists $k\in \Bbb{N}$\thinspace \thinspace $\setminus
\,\,\{1\}$ such that the $k$-th mixed cumulant $v_{0}$-compressed cumulants
of $x$ and $y$ does not vanish. Then since

\strut

$\ \ \ \ \ \ FP_{*}(G:x)$ $=$ $FP_{*}(G:a)\cap loop_{v_{0}}(G:a)$

and

$\ \ \ \ \ \ FP_{*}(G:y)$ $=$ $FP_{*}(G:b)\cap loop_{v_{0}}(G:b),$

\strut

The $k$-th mixed cumulant of $a$ and $b$ does not vanish. This contradict
our assumption.
\end{proof}

\strut

\strut \strut

\strut

\section{Diagonal Compressed Graph $W^{*}$-Probability}

\strut

\strut

Throughout this chapter, we will let $G$ be a countable directed graph and $%
\mathbb{F}^{+}(G)$, the free semigroupoid of $G$ and let $\left(
W^{*}(G),E\right) $ be the graph $W^{*}$-probability space over the diagonal
subalgebra $D_{G}.$ In this chapter, we will consider the diagonal
compression of a $D_{G}$-valued random variable $a\in \left(
W^{*}(G),E\right) ,$ for the given \textbf{finite} subset of the vertex set $%
V(G)$ of the graph $G.$ Let $V=\{v_{1},...,v_{N}\}$ be a subset of the
vertex set $V(G)$ of the graph $G.$ The diagonal compressed random variable
of $a\in \left( W^{*}(G),E\right) $ by $V$ is defined by

\strut

\begin{center}
$L_{v_{1}}aL_{v_{1}}+...+L_{v_{N}}aL_{v_{N}}$
\end{center}

\strut

in $\left( W^{\ast }(G),E\right) ,$ as a new $D_{G}$-valued random variable
in $\left( W^{\ast }(G),E\right) .$ Notice that each $L_{v_{j}}aL_{v_{j}}$
is the $v_{j}$-compressed random variable of $a,$ for \ $j=1,...,N,$ and it
can be regarded as a random variable in $(W_{v_{j}}^{\ast }(G),E_{v_{j}}),$
the $W^{\ast }$-probability space (over $\mathbb{C}$). In this chapter, we
will regard them as compressed $D_{G}$-valued random variables. By regarding
all $L_{v_{j}}aL_{v_{j}}$ as $D_{G}$-valued random variables in $\left(
W^{\ast }(G),E\right) ,$ the diagonal compressed random variable of $a$ by $%
V=\{v_{1},...,v_{N}\}\subset V(G)$ is also a $D_{G}$-valued random variable
in $\left( W^{\ast }(G),E\right) .$ Of course, in the subset $V$ of $V(G),$
the vertices satisfy that

\strut

\begin{center}
$v_{i}\neq v_{j},$ whenever $i\neq j$ in $\{1,...,N\}.$
\end{center}

\strut

In this chapter, we will observe the amalgamated ($D_{G}$-valued) free
probability information of such diagonal compressed random variables in $%
\left( W^{*}(G),E\right) .$

\strut

\begin{definition}
Let $V=\{v_{1},...,v_{N}\}\subset V(G).$ Define the diagonal compression by $%
V,$

\strut

$\ \ \ \ \ \ \ P_{V}:W^{*}(G)\rightarrow
\sum_{j=1}^{N}L_{v_{j}}W^{*}(G)L_{v_{j}}\subset W^{*}(G)$

by

$\ \ \ \ \ \ \ P_{V}(a)=\sum_{j=1}^{N}L_{v_{j}}aL_{v_{j}},$ \ for all \ $%
a\in W^{*}(G).$

\strut

We say that the $D_{G}$-valued random variable $P_{V}(a)$ is the diagonal
compressed random variable of $a$ by $V.$
\end{definition}

\strut

\strut

\strut

\subsection{Diagonal Compressed Moments and Cumulants}

\strut

\strut

Throughout this section, fix $N\geq 2$ in $\mathbb{N}$ and let $%
V=\{v_{1},...,v_{N}\}$ be the fixed subset of the vertex set $V(G).$ Let $%
a\in \left( W^{\ast }(G),E\right) $ be an arbitrary $D_{G}$-valued random
variable. Then, we can get the diagonal compressed random variable of $a$ by
the given set $V,$ $P_{V}(a)\in \left( W^{\ast }(G),E\right) .$ By the very
definition, we have that

\strut

\strut (4.1)

\begin{center}
$
\begin{array}{ll}
P_{V}(a) & =L_{v_{1}}aL_{v_{1}}+...+L_{v_{N}}aL_{v_{N}} \\ 
&  \\ 
& =\sum_{j=1}^{N}\,\underset{w=v_{j}av_{j}\in \mathbb{F}^{+}(G:a),\,u_{w}\in
\{1,*\}\,}{\sum }p_{w}L_{w}^{u_{w}},
\end{array}
$
\end{center}

\strut \strut \strut

where

\strut

\begin{center}
$a=\underset{w\in \mathbb{F}^{+}(G:a),\,u_{w}\in \{1,*\}}{\sum }%
p_{w}L_{w}^{u_{w}}\in \left( W^{*}(G),E\right) .$
\end{center}

\strut

If $V\subseteq V(G:a),$ then we have that

\strut

\begin{center}
$P_{V}(a)=\sum_{j=1}^{N}\left( [p_{v_{j}}L_{v_{j}}]+\underset{w\in
loop_{v_{j}}(G:a),\,u_{w}\in \{1,\ast \}}{\sum }p_{w}L_{w}^{u_{w}}\right) .$
\end{center}

\strut

But it is possible that $V\nsubseteqq V(G:a)$ and then $%
\sum_{j=1}^{N}[p_{v_{j}}L_{v_{j}}]=0_{D_{G}}.$ Choose $(i,j)\in
\{1,...,N\}^{2}$ such that $i\neq j$ and assume that there is at least one
loop $l=v_{i}lv_{i}\in loop_{v_{i}}(G:a)$ containing $v_{j}$ (where, $%
v_{i}\neq v_{j}\in V$). i.e.e, $l=v_{i}lv_{i}=v_{j}lv_{j}.$ Then

\strut

\begin{center}
$loop_{v_{i}}(G:a)\cap loop_{v_{j}}(G:a)\neq \emptyset .$
\end{center}

\strut

So, we need to be careful the intersections of $loop_{v_{j}}(G:a)$ are empty
or not.

\strut

\begin{lemma}
Let $v_{i}\neq v_{j}\in V(G)$ and assume that

\strut

\ $\ \ \ \ \ \ \ loop_{v_{i}}(G:a)\cap loop_{v_{j}}(G:a)=\emptyset .$

\strut

Then the vertex compressed random variables $L_{v_{i}}aL_{v_{i}}$ and $%
L_{v_{j}}aL_{v_{j}}$ of $a$ satisfy that

\strut

$\ \ \ \ \ \ \left( L_{v_{i}}aL_{v_{i}}\right) ^{m}\left(
L_{v_{j}}aL_{v_{j}}\right) ^{n}=0_{D_{G}},$ for all $m,n\in \mathbb{N}.$
\end{lemma}

\strut

\begin{proof}
Assume that $v_{i}\neq v_{j}$ in $V$ \ (i.e.e $i\neq j$) and suppose that

\strut

\ $\ \ \ \ \ \ \ loop_{v_{i}}(G:a)\cap loop_{v_{j}}(G:a)=\emptyset .$

\strut

Then we can easily conclude that there is no loop $l=v_{i}lv_{i}\in
loop(G:a) $ containing $v_{j}$ (or equivalently, there is no loop $%
l=v_{j}lv_{j}\in loop(G:a)$ containing $v_{i}$). Let $m=1=n.$ Then

\strut

\ \ $\ \ \ \ \ \ \left( L_{v_{i}}aL_{v_{i}}\right)
(L_{v_{j}}aL_{v_{j}})=L_{v_{i}}a\left( L_{v_{i}}L_{v_{j}}\right)
aL_{v_{j}}=0_{D_{G}}.$

\strut \strut

Similarly, if $m=n,$ then

\strut

\ \ $\ \ \ \ \ \ 
\begin{array}{ll}
(L_{v_{i}}aL_{v_{i}})^{m}(L_{v_{j}}aL_{v_{j}})^{m} & =\left(
L_{v_{i}}aL_{v_{i}}L_{v_{j}}aL_{v_{j}}\right) ^{m} \\ 
& =\left( L_{v_{i}}a(L_{v_{i}}L_{v_{j}})aL_{v_{j}}\right) ^{m}=0_{D_{G}}.
\end{array}
$

\strut

Now, let $m>n.$ Then

\strut

$\ \ \ \ \ \ 
\begin{array}{ll}
(L_{v_{i}}aL_{v_{i}})^{m}(L_{v_{j}}aL_{v_{j}})^{n} & 
=(L_{v_{i}}aL_{v_{i}})^{m-n}(L_{v_{i}}aL_{v_{i}})^{n}(L_{v_{j}}aL_{v_{j}})^{n}
\\ 
& =(L_{v_{i}}aL_{v_{i}})^{m-n}\left(
L_{v_{i}}a(L_{v_{i}}L_{v_{j}})aL_{v_{j}}\right) ^{n} \\ 
& =0_{D_{G}}.
\end{array}
$

\strut

Similarly, if $m<n,$ then $%
(L_{v_{i}}aL_{v_{i}})^{m}(L_{v_{j}}aL_{v_{j}})^{n}=0_{D_{G}}.$
\end{proof}

\strut

\begin{lemma}
Let $a\in \left( W^{*}(G),E\right) $ be a $D_{G}$-valued random variable and
let $v_{i}\neq v_{j}$ in $V$ and assume that

\strut

$\ \ \ \ loop_{v_{i}}^{v_{j}}(G:a)\overset{def}{=}loop_{v_{i}}(G:a)\cap
loop_{v_{j}}(G:a)\neq \emptyset .$

\strut

Then the vertex compressed random variables $L_{v_{i}}aL_{v_{i}}$ and $%
L_{v_{j}}aL_{v_{j}}$ satisfy that

\strut

$\ \ \ \left( L_{v_{i}}aL_{v_{i}}\right) ^{m}\left(
L_{v_{j}}aL_{v_{j}}\right) ^{n}$

\strut

$\ \ \ \ \ \ \ \ \ \ \ =\underset{(w_{1},...,w_{m+n})\in
loop_{v_{i}}^{v_{j}}(G:a)^{m+n}}{\sum }\left( \Pi
_{k=1}^{m+n}p_{w_{k}}\right)
L_{w_{1}}^{u_{w_{1}}}...L_{w_{m+n}}^{u_{w_{m+n}}}.$
\end{lemma}

\strut

\begin{proof}
By the assumption that

\strut

$\ \ \ \ \ \ \ \ \ loop_{v_{i}}(G:a)\cap loop_{v_{j}}(G:a)\neq \emptyset ,$

\strut

we can define the subset $loop_{v_{i}}^{v_{j}}(G:a)$ of $loop(G:a)\subset
FP(G:a)$ by

\strut

$\ \ \ 
\begin{array}{ll}
loop_{v_{i}}^{v_{j}}(G:a) & \overset{def}{=}loop_{v_{i}}(G:a)\cap
loop_{v_{j}}(G:a) \\ 
& =\{l\in loop(G:a):l=v_{i}lv_{i}\,=v_{j}lv_{j}\}.
\end{array}
$

\strut

By the very definition, we have that

\strut

$\ \ \ \ \ \ \ \ \ \ \ \ \
loop_{v_{i}}^{v_{j}}(G:a)=loop_{v_{j}}^{v_{i}}(G:a)$

\strut

in $loop(G:a).$ By the observation in Section 1.3, we have that

\strut

$\ \ \left( L_{v_{i}}aL_{v_{i}}\right) ^{m}=\underset{(w_{1},...,w_{m})\in
loop_{v_{i}}(G:a)^{m},\,u_{w_{k}}\in \{1,*\}}{\sum }\left( \Pi
_{k=1}^{n}p_{w_{k}}\right) L_{w_{1}}^{u_{w_{1}}}...L_{w_{m}}^{u_{w_{m}}}$

\strut

and

\strut

$\ \ \left( L_{v_{j}}aL_{v_{j}}\right) ^{n}=\underset{(w_{1}^{\prime
},...,w_{n}^{\prime })\in loop_{v_{j}}(G:a)^{n},\,u_{w_{p}^{\prime }}\in
\{1,*\}}{\sum }\left( \Pi _{p=1}^{n}p_{w_{p}^{\prime }}\right)
L_{w_{1}^{\prime }}^{u_{w_{1}^{\prime }}}...L_{w_{n}^{\prime
}}^{u_{w_{n}^{\prime }}}.$

\strut

\strut

Thus

\strut

$\ (L_{v_{i}}aL_{v_{i}})^{m}(L_{v_{j}}aL_{v_{j}})^{n}$

\strut

$\ =\left( \underset{(w_{1},...,w_{m})\in
loop_{v_{i}}(G:a)^{m},,\,u_{w_{k}}\in \{1,*\}}{\sum }\left( \Pi
_{k=1}^{n}p_{w_{k}}\right)
L_{w_{1}}^{u_{w_{1}}}...L_{w_{m}}^{u_{w_{m}}}\right) $

\strut

$\ \ \ \ \ \ \ \left( \underset{(w_{1}^{\prime },...,w_{n}^{\prime })\in
loop_{v_{j}}(G:a)^{n},\,u_{w_{p}^{\prime }}\in \{1,*\}}{\sum }\left( \Pi
_{p=1}^{n}p_{w_{p}^{\prime }}\right) L_{w_{1}^{\prime }}^{u_{w_{1}^{\prime
}}}...L_{w_{n}^{\prime }}^{u_{w_{n}^{\prime }}}\right) $

\strut

$\ =\underset{(w_{1},...,w_{m},w_{1}^{\prime },...,w_{n}^{\prime })\in
loop_{v_{i}}^{v_{j}}(G:a)^{m+n}}{\sum }\left( \Pi _{k=1}^{m}p_{w_{k}}\right)
\left( \Pi _{p=1}^{n}p_{w_{p}^{\prime }}\right) $

\strut

$\ \ \ \ \ \ \ \ \ \ \ \ \ \ \ \ \ \ \ \ \ \ \ \left(
L_{w_{1}}^{u_{w_{1}}}...L_{w_{m}}^{u_{w_{m}}}L_{w_{1}^{\prime
}}^{u_{w_{1}^{\prime }}}...L_{w_{n}^{\prime }}^{u_{w_{n}^{\prime }}}\right)
. $

\strut
\end{proof}

\strut

The above lemma says that, in general, if $loop_{v_{i}}^{v_{j}}(G:a)\neq
\emptyset ,$ then

\strut

\begin{center}
$(L_{v_{i}}aL_{v_{i}})(L_{v_{j}}aL_{v_{j}})=\underset{(w,w^{\prime })\in
loop_{v_{i}}^{v_{j}}(G:a)^{2}}{\sum }p_{w}p_{w^{\prime
}}L_{w}^{u_{w}}L_{w^{\prime }}^{u_{w^{\prime }}},$
\end{center}

\strut

where $u_{w},u_{w^{\prime }}\in \{1,\ast \}$ and $a=\underset{w\in
FP(G:a),\,u_{w}\in \{1,\ast \}}{\sum }p_{w}L_{w}^{u_{w}}\in \left( W^{\ast
}(G),E\right) .$

\strut

Now, let $d_{k}=\underset{v^{(k)}\in V(G:d_{k})}{\sum }%
q_{v^{(k)}}L_{v^{(k)}}\in D_{G},$ for $k\in \mathbb{N}.$ Then

\strut

$\ \ d_{1}\left( P_{V}(a)\right) d_{2}\left( P_{V}(a)\right) ...d_{n}\left(
P_{V}(a)\right) $

\strut

$\ \ \ \ =\underset{(v^{(1)},...,v^{(n)})\in \Pi _{k=1}^{n}V(G:d_{k})}{\sum }%
\left( \Pi _{k=1}^{n}~q_{v^{(k)}}\right) \left( L_{v^{(1)}}\left(
P_{V}(a)\right) ~...~L_{v^{(n)}}(P_{V}(a))\right) $

\strut

$\ \ \ \ =\underset{(v^{(1)},...,v^{(n)})\in \Pi _{k=1}^{n}V(G:d_{k})}{\sum }%
\,\underset{(w_{1},...,w_{n})\in
loop_{v_{j}}(G:P_{V}(a))^{n},\,\,u_{w_{j}}\in \{1,*\}}{\sum }$

\strut

$\ \ \ \ \ \ \ \ \ \ \ \ \ \ \ \ \ \ \ \ \ \ \ \ \ \ \ \ \ \ \ \ \ \ \ \ \ \
\ \left( \Pi _{j=1}^{n}q_{v^{(j)}}\right) \left( \Pi
_{j=1}^{n}p_{w_{j}}\right) $

\strut

$\ \ \ \ \ \ \ \ \ \ \ \ \ \ \ \ \ \ \ \ \ \ \ \ \ \ \ \ \ \ \ \ \ \ \ \ \ \
\ \ \ \ \ \ \ \ \ \ \ \ \ \ \ \ \ \left( \Pi _{j=1}^{n}\delta
_{v^{(j)},\,x_{j}}\right) L_{w_{1}}^{u_{w_{1}}}...L_{w_{n}}^{u_{w_{n}}}.$

\strut

We have that

$\strut $

\begin{center}
$\mathbb{F}^{+}(G:P_{V}(a))=\left( \cup _{k=1}^{n}loop_{v_{k}}(G:a)\right)
\cup \left( V\cap V(G:a)\right) ,$
\end{center}

\strut

where

\begin{center}
$P_{V}(a)=\underset{w\in \mathbb{F}^{+}(G:P_{V}(a)),\,u_{w}\in \{1,*\}}{\sum 
}p_{w}L_{w}^{u_{w}}\in \left( W^{*}(G),E\right) $
\end{center}

and

\begin{center}
$a=\underset{w\in \mathbb{F}^{+}(G:a),\,u_{w}\in \{1,*\}}{\sum }%
p_{w}L_{w}^{u_{w}}\in \left( W^{*}(G),E\right) .$
\end{center}

\strut

Notice that the coefficients $p_{w}=\,<\xi _{w},a\xi _{w}>$'s are not
changed, because $\mathbb{F}^{+}(G:P_{V}(a))\subset \mathbb{F}^{+}(G:a).$
Now, we have all information to get the $D_{G}$-valued moments of the
diagonal compressed random variable of $a$ by $V$ ;

\strut

\begin{theorem}
Let $a\in \left( W^{*}(G),E\right) $ be a $D_{G}$-valued random variable and
let $V=\{v_{1},...,v_{N}\}$ be the fixed finite subset of the vertex set $%
V(G).$ The diagonal compressed random variable of $a$ by $V,$ $P_{V}(a)$ has
the $n$-th moment

\strut

$E\left( d_{1}P_{V}(a)d_{2}P_{V}(a)...d_{n}P_{V}(a)\right) $

\strut

$\ =\underset{\pi \in NC(n)}{\sum }\,\underset{(v^{(1)},...,v^{(n)})\in \Pi
_{k=1}^{n}V(G:d_{k})}{\sum }\,$

\strut

$\ \ \ \ \ \ \ \underset{(w_{1},...,w_{n})\in \left( \left( \cup
_{k=1}^{N}loop_{v_{k}}(G:a)\right) \cup \left( V\cap V(G:a)\right) \right)
^{n},\,w_{j}=x_{j}w_{j}x_{j},\,\,u_{w_{j}}\in \{1,*\}}{\sum }$

$\strut $

$\ \ \ \ \ \ \ \ \ \ \ \ \ \ \ \ \ \ \ \ \ \left( \Pi
_{j=1}^{n}q_{v^{(j)}}\right) \left( \Pi _{j=1}^{n}p_{w_{j}}\right) \left(
\Pi _{j=1}^{n}\delta _{v^{(j)},\,x_{j}}\right) $

$\strut $

$\ \ \ \ \ \ \ \ \ \ \ \ \ \ \ \ \ \ \ \ \ \ \ \ \ \ \ \ \ \ \ \ \ \ \ \ \ \
\ \ \ E\left( L_{w_{1}}^{u_{w_{1}}}...L_{w_{n}}^{u_{w_{n}}}\right) ,$

\strut

for all $n\in \mathbb{N},$ where $d_{k}=\underset{v^{(k)}\in V(G:d_{k})}{%
\sum }q_{v^{(k)}}L_{v^{(k)}}\in D_{G}$ are arbitrary, for $k=1,...,n.$
\end{theorem}

\strut

\begin{proof}
Fix $n\in \mathbb{N}$ and let $d_{k}=\underset{v^{(k)}\in V(G:d_{k})}{\sum }%
q_{v^{(k)}}L_{v^{(k)}}\in D_{G}$ are arbitrary, for $k=1,...,n.$ By using
the same notations in Section 1.3, we have that

\strut

$\ E\left( d_{1}P_{V}(a)...d_{n}P_{V}(a)\right) $

\strut

$\ \ =\underset{(v^{(1)},...,v^{(n)})\in \Pi _{k=1}^{n}V(G:d_{k})}{\sum }\,%
\underset{(w_{1},...,w_{n})\in \mathbb{F}^{+}(G:P_{V}(a))^{n},%
\,w_{j}=x_{j}w_{j}x_{j},\,\,u_{w_{j}}\in \{1,\ast \}}{\sum }$

\strut

$\ \ \ \ \ \ \ \ \ \ \ \ \ \ \ \ \ \ \ \ \ \ \ \ \ \ \ \ \ \ \ \ \ \ \ \ \ \
\left( \Pi _{j=1}^{n}q_{v^{(j)}}\right) \left( \Pi
_{j=1}^{n}p_{w_{j}}\right) $

\strut

$\ \ \ \ \ \ \ \ \ \ \ \ \ \ \ \ \ \ \ \ \ \ \ \ \ \ \ \ \ \ \ \ \ \ \ \ \ \
\left( \Pi _{j=1}^{n}\delta _{v^{(j)},\,x_{j}}\right) E\left(
L_{w_{1}}^{u_{w_{1}}}...L_{w_{n}}^{u_{w_{n}}}\right) $

\strut

$\ \ =\underset{\pi \in NC(n)}{\sum }\,\underset{(v^{(1)},...,v^{(n)})\in
\Pi _{k=1}^{n}V(G:d_{k})}{\sum }\,\underset{(w_{1},...,w_{n})\in \mathbb{F}%
^{+}(G:P_{V}(a))^{n},\,w_{j}=x_{j}w_{j}x_{j},\,\,u_{w_{j}}\in \{1,\ast \}}{%
\sum }$

\strut \strut

$\ \ \ \ \ \ \ \ \ \ \ \ \ \ \ \ \ \ \ \ \ \ \ \ \ \ \ \ \ \ \ \ \ \ \ \
\left( \Pi _{j=1}^{n}q_{v^{(j)}}\right) \left( \Pi
_{j=1}^{n}p_{w_{j}}\right) $

\strut \strut

$\ \ \ \ \ \ \ \ \ \ \ \ \ \ \ \ \ \ \ \ \ \ \ \ \ \ \ \ \ \ \ \ \ \ \ \
\left( \Pi _{j=1}^{n}\delta _{v^{(j)},\,x_{j}}\right) \,\,\,\Pr oj\left(
L_{w_{1}}^{u_{w_{1}}}...L_{w_{n}}^{u_{w_{n}}}\right) ,$

\strut \strut \strut

where

\strut

$\ \ \ \ \ \ \ \mathbb{F}^{+}(G:P_{V}(a))=\left( \cup
_{k=1}^{N}loop_{v_{k}}(G:a)\right) \cup \left( V\cap V(G:a)\right) .$

\strut
\end{proof}

\strut

By the M\"{o}bius inversion, we can get the $n$-th cumulants of $P_{V}(a)$ ;

\strut

\begin{theorem}
Let $a\in \left( W^{*}(G),E\right) $ be a $D_{G}$-valued random variable and
let $V=\{v_{1},...,v_{N}\}$ be the fixed subset of $V(G).$ Then the diagonal
compressed random variable of $a$ by $V,$ $P_{V}(a)$ has the $n$-th
cumulants are

\strut

$\ \ \ \ \ \ \ k_{1}\left( d_{1}P_{V}(a)\right) =\underset{v\in V\cap
(V(G:d_{1})\cap V(G:a))}{\sum }\left( q_{v}p_{v}\right) L_{v}$

\strut and

\strut \strut

$\ \ k_{n}\left( \underset{n-times}{\underbrace{%
d_{1}P_{V}(a),....,d_{n}P_{V}(a)}}\right) $

\strut

$\ \ \ \ \ =\underset{(v^{(1)},...,v^{(n)})\in \Pi _{k=1}^{n}V(G:d_{k})}{%
\sum }\,\left( \Pi _{j=1}^{n}q_{v^{(j)}}\right) $

\strut \strut

$\ \ \ \ \ \ \ \underset{(w_{1},...,w_{n})\in \left( \left( \cup
_{k=1}^{n}loop_{v_{k}}(G:a)\right) \cup \left( V\cap V(G:a)\right) \right)
^{n},~\,\,u_{w_{j}}\in \{1,*\}}{\sum }\left( \Pi _{j=1}^{n}p_{w_{j}}\right) $

\strut

$\ \ \ \ \ \ \ \ \ \ \ \ \ \ \ \ \ \ \ \left( \Pi _{k=1}^{n}\delta
_{v^{(k)},~x_{k}}\right) \,\,\,\ \mu
_{w_{1},...,w_{n}}^{u_{w_{1}},...,u_{w_{n}}}E\left(
L_{w_{1}}^{u_{w_{1}}}...L_{w_{n}}^{u_{w_{n}}}\right) ,$

\strut

for all $n>1$ in $\mathbb{N},$ where $d_{k}=\underset{v^{(k)}\in V(G:d_{k})}{%
\sum }q_{v^{(k)}}L_{v^{(k)}}\in D_{G}$ are arbitrary for $k=1,...,n.$
\end{theorem}

\strut

\begin{proof}
Let $n=1.$ Then

\strut

$\ \ \ k_{1}\left( d_{1}P_{V}(a)\right) =d_{1}E\left( P_{V}(a)\right) $

\strut

$\ \ \ \ \ \ \ =d_{1}E\left( \underset{v\in V\cap V(G:a)}{\sum }p_{v}L_{v}+%
\underset{w\in \cup _{k=1}^{N}loop_{v_{k}}(G:a),\,u_{w}\in \{1,*\}}{\sum }%
p_{w}L_{w}^{u_{w}}\right) $

\strut

$\ \ \ \ \ \ \ =\left( \underset{v^{(1)}\in V(G:d_{1})}{\sum }%
q_{v^{(1)}}L_{v^{(1)}}\right) \left( \underset{v\in V\cap V(G:a)}{\sum }%
p_{v}L_{v}\right) $

\strut

$\ \ \ \ \ \ \ =\underset{v\in V\cap (V(G:d_{1})\cap V(G:a))}{\sum }\left(
q_{v}p_{v}\right) L_{v}.$

\strut

Now, fix $n>1$ in $\mathbb{N}.$ By the M\"{o}bius inversion and by Section
1.3, we have that

\strut

$\ k_{\pi }\left( d_{1}P_{V}(a),...,d_{n}P_{V}(a)\right) $

\strut

$\ \ =\underset{(v^{(1)},...,v^{(n)})\in \Pi _{k=1}^{n}V(G:d_{k})}{\sum }\,%
\underset{(w_{1},...,w_{n})\in \mathbb{F}^{+}(G:P_{V}(a))^{n},%
\,w_{j}=x_{j}w_{j}x_{j},\,\,u_{w_{j}}\in \{1,*\}}{\sum }$

\strut

$\ \ \ \ \ \ \ \ \ \ \ \ \ \ \ \ \ \ \ \ \ \ \ \left( \Pi
_{j=1}^{n}q_{v^{(j)}}\right) \left( \Pi _{j=1}^{n}p_{w_{j}}\right) \left(
\Pi _{j=1}^{n}\delta _{v^{(j)},\,x_{j}}\right) \,\,\,$

\strut

$\ \ \ \ \ \ \ \ \ \ \ \ \ \ \ \ \ \ \ \ \ \ \ \ \ \ \ \ \ \ \ \ \ \ \ \ \ \
\ \ \ \ \ \mu _{w_{1},...,w_{n}}^{u_{w_{1}},...,u_{w_{n}}}E\left(
L_{w_{1}}^{u_{w_{1}}}...L_{w_{n}}^{u_{w_{n}}}\right) ,$

\strut

where

\strut

$\ \ \mathbb{F}^{+}(G:P_{V}(a))=\left( \cup
_{k=1}^{n}loop_{v_{k}}(G:a)\right) \cup \left( V\cap V(G:a)\right) .$

\strut
\end{proof}

\strut

Notice that in the set $\cup _{k=1}^{N}loop_{v_{k}}(G:a),$ it is possible
that there are subsets $loop_{v_{i}}^{v_{j}}(G:a),$ for $i\neq j$ in $%
\{1,...,N\}.$

\strut

Now, assume that the fixed subset $V=\{v_{1},...,v_{N}\}\subset V(G)$
satisfies that, for any choice $\left( v_{i},v_{j}\right) \in V^{2}$ with $%
i\neq j,$ there is no loop $l=v_{i}lv_{i}$ containing $v_{j}$ in the graph $%
G.$ Notice that, in this case, $W_{v_{i}}^{*}(G)$ and $W_{v_{j}}^{*}(G)$ are
free over $D_{G},$ in $\left( W^{*}(G),E\right) .$ Therefore, we can
conclude that $W_{v_{1}}^{*}(G)$,..., $W_{v_{N}}^{*}(G)$ are mutually free
over $D_{G},$ in $\left( W^{*}(G),E\right) ,$ for the given subset $%
V=\{v_{1},...,v_{N}\}$ in $V(G).$ Also, in this case, the diagonal
compression of a $D_{G}$-valued random variable $a\in \left(
W^{*}(G),E\right) ,$ by $V,$ \ $P_{V}(a)$ can be regarded as the free sum of
vertex compressed random variables of $a$ by $v_{j}$'s, \ $j=1,...,N.$ i.e.e
the diagonal compressed random variable of $a$ by $V$

\strut

\begin{center}
$P_{V}(a)=L_{v_{1}}aL_{v_{1}}+...+L_{v_{N}}aL_{v_{N}}$
\end{center}

\strut

and the vertex compressed random variables $%
L_{v_{1}}aL_{v_{1}},...,L_{v_{N}}aL_{v_{N}}$ are mutually free over $D_{G},$
in $\left( W^{*}(G),E\right) .$ Thus $P_{V}(a)$ is the sum of $D_{G}$-free
random variables $L_{v_{1}}aL_{v_{1}},...,L_{v_{N}}aL_{v_{N}}.$

\strut

\begin{lemma}
Let $a\in \left( W^{*}(G),E\right) $ be a $D_{G}$-valued random variable and
let $a_{v_{0}}=L_{v_{0}}aL_{v_{0}}$ be the vertex compressed random variable
of $a$ by the fixed vertex $v_{0}\in V(G).$ Then $da_{v_{0}}=a_{v_{0}}d,$
for all $d\in D_{G}.$
\end{lemma}

\strut

\begin{proof}
Let $a_{v_{0}}=L_{v_{0}}aL_{v_{0}}\in \left( W^{*}(G),E\right) $ be the $v$%
-compressed random variable and let $d=\underset{v\in V(G:d)}{\sum }%
q_{v}L_{v}\in D_{G}$ be arbitrary. We can express $a_{v_{0}}$ by

\strut

$\ \ \ \ \ \ \ \ \ \ \ a_{v_{0}}=[p_{v_{0}}L_{v_{0}}]+\underset{w\in
loop_{v_{0}}(G:a),\,u_{w}\in \{1,\ast \}}{\sum }p_{w}L_{w}^{u_{w}},$

\strut

where the $D_{G}$-valued random variable $a$ has the form

\strut

$\ \ \ \ \ \ \ \ \ \ \ a=\underset{w\in \mathbb{F}^{+}(G:a),\,u_{w}\in
\{1,\ast \}}{\sum }p_{w}L_{w}^{u_{w}}\in \left( W^{\ast }(G),E\right) $

\strut

and where

\strut

$\ \ \ \ \ \ \ \ \ \ \ [p_{v_{0}}L_{v_{0}}]=\left\{ 
\begin{array}{lll}
p_{v_{0}}L_{v_{0}} &  & \text{if }v_{0}\in V(G:a) \\ 
&  &  \\ 
0_{D_{G}} &  & \text{otherwise.}
\end{array}
\right. $

\strut

Suppose that $v_{0}\notin V(G:d).$ Then $da_{v_{0}}=0_{D_{G}}=a_{v_{0}}d.$
Now assume that $v_{0}\in V(G:d).$ Then

\strut

$\ \ \ da_{v_{0}}=\left( q_{v_{0}}L_{v_{0}}+D\right)
a_{v_{0}}=q_{v_{0}}L_{v_{0}}a_{v_{0}}+Da_{v_{0}}$

$\strut $

$\ \ \ \ \ \ \ \ \ =q_{v_{0}}a_{v_{0}}+0_{D_{G}}$

$\strut $

$\ \ \ \ \ \ \ \ \ =q_{v_{0}}a_{v_{0}}+a_{v_{0}}D=a_{v_{0}}\left(
q_{v_{0}}L_{v_{0}}+D\right) $

$\strut $

$\ \ \ \ \ \ \ \ \ =a_{v_{0}}d,$

\strut

for all $d=q_{v_{0}}L_{v_{0}}+D\in D_{G}$ having its summand $%
q_{v_{0}}L_{v_{0}},$ since $[p_{v_{0}}L_{v_{0}}]$ commutes with $d$ and $%
\underset{w\in loop_{v_{0}}(G:a),\,u_{w}\in \{1,\ast \}}{\sum }%
p_{w}L_{w}^{u_{w}}$ also commutes with $d.$
\end{proof}

\strut \strut

By the previous lemma, we can conclude that ;

\strut

\begin{proposition}
Let $a\in \left( W^{*}(G),E\right) $ be a $D_{G}$-valued random variable and
let $V=\{v_{1},...,v_{N}\}$ be a finite subset of the vertex set $V(G).$
Assume that

\strut

$\ \ \ \ \ \ \ \ \ \ \ \ \ \ \ \ \ \ \ \ \ \ \ \ \
loop_{v_{i}}^{v_{j}}(G:a)=\emptyset ,$

\strut

for any $(v_{i},v_{j})\in V^{2}$ such that $i\neq j$ in $\{1,...,N\}.$ Then
the diagonal compressed random variable of $a$ by $V,$ \ $P_{V}(a)$ has $n$%
-th cumulants

\strut

$\ \ \ \ k_{n}\left( d_{1}P_{V}(a),...,d_{n}P_{V}(a)\right) $

\strut

$\ \ \ \ \ \ =\left\{ 
\begin{array}{cc}
\sum_{j=1}^{N}\,\,(d_{1}d_{2})\left( \underset{l\in loop_{*}^{v_{j}}(G:a)}{%
\sum }2(p_{l}p_{l^{t}})^{2}\cdot L_{v_{j}}\right) & \text{if \ }n=2 \\ 
&  \\ 
0_{D_{G}} & \text{otherwise,}
\end{array}
\right. ,$

\strut

for all $nN,$ where $d_{k}=\underset{v^{(k)}\in V(G:d_{k})}{\sum }%
q_{v^{(k)}}L_{v^{(k)}}\in D_{G}$ are arbitrary, \ $j=1,...,n.$
\end{proposition}

\strut

\begin{proof}
Fix $n\in \mathbb{N}$ and denote $L_{v_{j}}aL_{v_{j}}$ by $a_{j},$ for \ $%
j=1,...,N.$ Then we have that

\strut

$\ \ \ \ \ \ \ \ \ \ \ \ \ FP_{\ast }(G:a_{j})=loop_{\ast }(G:a_{j})$

and

$\ \ \ \ \ \ \ \ \ \ \ \ \ FP_{\ast }^{c}(G:a_{j})=loop_{\ast
}^{c}(G:a_{j}), $

\strut

for all \ $j=1,...,N.$ By the assumption that $loop_{v_{i}}^{v_{j}}(G)=%
\varnothing ,$ for all $i\neq j\in \{1,...,N\},$ we can conclude that

\strut

$\ \ \ \ \ \ \ \ \ \ \ FP_{\ast }(G:a_{i})\cap FP_{\ast
}(G:a_{j})=\varnothing $

and

$\ \ \ \ \ \ \ \ \ \ \ FP_{\ast }^{c}(G:a_{i})\cap FP_{\ast
}^{c}(G:a_{j})=\varnothing ,$

\strut

for any $i\neq j\in \{1,...,N\}.$ So, $W_{*}^{\{a_{i},a_{j}\}}=\varnothing ,$
by the above second intersection.. This shows that the $D_{G}$-valued random
variables $a_{1},...,a_{N}$ are free from each other over $D_{G},$ in $%
(W^{*}(G),E).$ Thus our diagonal compressed random variable $P_{V}(a)$ by $V$
is the $D_{G}$-free sum of $a_{1},...,a_{N}.$ Then, for arbitrary $%
d_{1},...,d_{n}\in D_{G},$

\strut \strut

$\ k_{n}\left( d_{1}P_{V}(a),...,d_{n}P_{V}(a)\right) $

\strut

$\ \ \ =k_{n}\left( d_{1}(a_{1}+...+a_{N}),...,d_{n}(a_{1}+...+a_{n})\right) 
$

\strut

$\ \ \ =\sum_{j=1}^{N}k_{n}\left( d_{1}a_{j},...,d_{n}a_{j}\right) $

\strut \strut

by the mutual $D_{G}$-freeness of $a_{1},...,a_{N}$ in $\left(
W^{*}(G),E\right) $

\strut

$\ \ \ =\sum_{j=1}^{N}\,\,(d_{1}...d_{n})\,\,k_{n}\left( \underset{n-times}{%
\underbrace{a_{j},.....,a_{j}}}\right) $

\strut

by the previous lemma and by the bimodule map property of the $D_{G}$-valued
cumulant

\strut

$\ \ \ =\sum_{j=1}^{N}\,\,(d_{1}...d_{n})~k_{n}\left( a_{j~(\ast
)},...,a_{j~(\ast )}\right) $

\strut

where $a_{j}=a_{j~d}+a_{j~(\ast )}+a_{j~(non-\ast )}$

\strut

$\ \ \ =\sum_{j=1}^{N}\,\,(d_{1}...d_{n})~\left( \underset{l\in loop_{\ast
}^{v_{j}}(G:a)}{\sum }k_{n}\left( p_{l}L_{l}+p_{l^{t}}L_{l}^{\ast
}~,...,~p_{l}L_{l}+p_{l^{t}}L_{l}^{\ast }\right) \right) $

\strut

by Section 3.2

\strut

$\ \ \ =\sum_{j=1}^{N}\,\,(d_{1}...d_{n})\underset{l\in loop_{\ast
}^{v_{j}}(G:a)}{\sum }\left( p_{l}p_{l^{t}}\right) ^{n}\cdot $

\strut

$\ \ \ \ \ \ \ \ \ \ \ \ \ \ \ \ \ \ \ \ \ \ \ \ \ \ \underset{%
(u_{1},...,u_{n})\in \{1,\ast \}^{n}}{\sum }k_{n}\left(
L_{l}^{u_{1}},...,L_{l}^{u_{n}}\right) $

\strut

$\ \ \ =\sum_{j=1}^{N}\,\,(d_{1}...d_{n})\underset{l\in loop_{\ast
}^{v_{j}}(G:a)}{\sum }~\left( p_{l}p_{l^{t}}\right) ^{n}\cdot $

$\strut $

$\ \ \ \ \ \ \ \ \ \ \ \ \ \ \ \ \ \ \ \ \ \underset{L\in LP_{n}^{*}}{\sum }%
\mu _{l,...,l}^{L(u_{1},...,u_{n})}E\left(
L_{l}^{u_{1}}...L_{l}^{u_{n}}\right) $

\strut

$\ \ \ =\sum_{j=1}^{N}\,\,(d_{1}...d_{n})\left( \underset{l\in loop_{\ast
}^{v_{j}}(G:a)}{\sum }~\left( p_{l}p_{l^{t}}\right) ^{n}\cdot \underset{L\in
LP_{n}^{\ast }}{\sum }\mu _{l,...,l}^{L(u_{1},...,u_{n})}L_{v_{j}}\right) $

\strut

since $E(L_{l}^{u_{1}}...L_{l}^{u_{n}})=L_{v_{j}}$, for the suitable $%
(u_{1},...,u_{n})\in \{1,*\}^{n}$

\strut

$\ \ \ =\left\{ 
\begin{array}{cc}
\sum_{j=1}^{2}\,\,(d_{1}d_{2})\left( \underset{l\in loop_{\ast }^{v_{j}}(G:a)%
}{\sum }2(p_{l}p_{l^{t}})^{2}\cdot L_{v_{j}}\right) & \text{if \ }n=2 \\ 
&  \\ 
0_{D_{G}} & \text{otherwise,}
\end{array}
\right. $

\strut

by the $D_{G}$-semicircularity of $a_{j}$'s \ ($j=1,...,N$).
\end{proof}

\strut \bigskip

\bigskip \strut

\subsection{Diagonal Compressed Freeness on $\left( W^{*}(G),E\right) $}

\strut

\strut

In this section, we will consider the $D_{G}$-freeness of two diagonal
compressed random variables. Likewise the vertex compressed case, we can get
that ;

\strut

\begin{proposition}
Let $V=\{v_{1},...,v_{N}\}$ be a finite subset of the vertex set $V(G).$ Let 
$a,b\in \left( W^{*}(G),E\right) $ be $D_{G}$-valued random variables. If $a$
and $b$ are free over $D_{G},$ in $\left( W^{*}(G),E\right) ,$ then the
corresponding diagonal compressed random variables $P_{V}(a)$ and $P_{V}(b)$
of $a$ and $b$ by $V$ are free over $D_{G},$ in $\left( W^{*}(G),E\right) .$ 
$\square $
\end{proposition}

\strut \strut \strut \strut \strut \strut

\strut Now, we will consider the two diagonal compressions $P_{V_{1}}$ and $%
P_{V_{2}}$ of $D_{G}$-valued random variables in the graph $W^{*}$%
-probability space $\left( W^{*}(G),E\right) $ over the diagonal subalgebra $%
D_{G}.$ The $D_{G}$-freeness of two $D_{G}$-valued random variables $%
P_{V_{1}}(a)$ and $P_{V_{2}}(a)$ is determined as follows ;

\strut

\begin{proposition}
Let $V_{1}=\{v_{1}^{(1)},...,v_{N_{1}}^{(1)}\}$ and $V_{2}=%
\{v_{1}^{(2)},...,v_{N_{2}}^{(2)}\}$ be finite subsets of the vertex set $%
V(G),$ where $N_{1},N_{2}\in \mathbb{N}.$ Let $a\in \left( W^{*}(G),E\right) 
$ be an arbitrary $D_{G}$-valued random variable and let $P_{V_{1}}(a)$ and $%
P_{V_{2}}(a)$ be the corresponding diagonal compressed random variables by $%
V_{1}$ and $V_{2},$ respectively. If $V_{1}$ and $V_{2}$ satisfy that

\strut

$\ \ \ \ \ \ \ \ \ \ \ \ \ \ \ \ \ \ \ \ \ \ \ \ \ \ \ V_{1}\cap
V_{2}=\emptyset $

\strut and

$\ \ \ \ \ \ \ \ \ \ \ loop_{v_{i}^{(1)}}(G:a)\cap
loop_{v_{j}^{(2)}}(G:a)=\emptyset ,$

\strut

for all choices $(i,j)\in \{1,...,N_{1}\}\times \{1,...,N_{2}\},$ then $%
P_{V_{1}}(a)$ and $P_{V_{2}}(a)$ are free over $D_{G},$ in $\left(
W^{*}(G),E\right) .$
\end{proposition}

\strut

\begin{proof}
\strut Let $a=\underset{w\in \Bbb{F}^{+}(G:a),\,u_{w}\in \{1,*\}}{\sum }%
p_{w}L_{w}^{u_{w}}$ be a $D_{G}$-valued random variable in $\left(
W^{*}(G),E\right) $ and let $P_{V_{1}}$ and $P_{V_{2}}$ be the diagonal
compressions by $V_{1}$ and $V_{2},$ respectively. Suppose that $V_{1}\cap
V_{2}=\emptyset $ and assume that

\strut

$\ \ \ \ \ \ \ \ \ \ \ loop_{v_{i}^{(1)}}(G:a)\cap
loop_{v_{j}^{(2)}}(G:a)=\emptyset ,$

\strut

for all pairs $(i,j)\in \{1,...,N_{1}\}\times \{1,...,N_{2}\}.$ Then $l_{1}$
and $l_{2}$ are diagonal-distinct, for all $(i,j)$ $\in $ $\{1,...,N_{1}\}$ $%
\times $ $\{1,...,N_{2}\}.$ Thus we have that loops $l$ and $l^{\prime }$
are diagram-distinct, for all $(l,l^{\prime })$ $\in $ $loop\left(
G:P_{V_{1}}(a)\right) $ $\times $ $loop\left( G:P_{V_{2}}(a)\right) .$ By
the definition of diagonal compression, we also have that

\strut

$\ \ \ loop^{c}\left( G:P_{V_{1}}(a)\right) \cap loop^{c}\left(
G:P_{V_{2}}(a)\right) =\emptyset \cap \emptyset =\emptyset .$

\strut

Therefore, two $D_{G}$-valued random variables $P_{V_{1}}(a)$ and $%
P_{V_{2}}(a)$ are free over $D_{G}$ in $\left( W^{*}(G),E\right) .$\strut
\end{proof}

\strut \strut

\strut

\strut \bigskip \bigskip

\begin{quote}
\textbf{Reference}

\strut

\strut

{\small [1] \ \ A. Nica, R-transform in Free Probability, IHP course note,
available at www.math.uwaterloo.ca/\symbol{126}anica.}

{\small [2]\strut \ \ \ A. Nica and R. Speicher, R-diagonal Pair-A Common
Approach to Haar Unitaries and Circular Elements, (1995), www
.mast.queensu.ca/\symbol{126}speicher.\strut }

{\small [3] \ }$\ ${\small B. Solel, You can see the arrows in a Quiver
Operator Algebras, (2000), preprint}

{\small \strut [4] \ \ A. Nica, D. Shlyakhtenko and R. Speicher, R-cyclic
Families of Matrices in Free Probability, J. of Funct Anal, 188 (2002),
227-271.}

{\small [5] \ \ D. Shlyakhtenko, Some Applications of Freeness with
Amalgamation, J. Reine Angew. Math, 500 (1998), 191-212.\strut }

{\small [6] \ \ D.Voiculescu, K. Dykemma and A. Nica, Free Random Variables,
CRM Monograph Series Vol 1 (1992).\strut }

{\small [7] \ \ D. Voiculescu, Operations on Certain Non-commuting
Operator-Valued Random Variables, Ast\'{e}risque, 232 (1995), 243-275.\strut 
}

{\small [8] \ \ D. Shlyakhtenko, A-Valued Semicircular Systems, J. of Funct
Anal, 166 (1999), 1-47.\strut }

{\small [9] \ \ D.W. Kribs and M.T. Jury, Ideal Structure in Free
Semigroupoid Algebras from Directed Graphs, preprint}

{\small [10]\ D.W. Kribs and S.C. Power, Free Semigroupoid Algebras, preprint%
}

{\small [11]\ I. Cho, Amalgamated Boxed Convolution and Amalgamated
R-transform Theory, (2002), preprint.}

{\small [12] I. Cho, The Tower of Amalgamated Noncommutative Probability
Spaces, (2002), Preprint.}

{\small [13]\ I. Cho, Compatibility of a Noncommutative Probability Space
and a Noncommutative Probability Space with Amalgamation, (2003), Preprint}

{\small [14] I. Cho, An Example of Scalar-Valued Moments, Under
Compatibility, (2003), Preprint.}

{\small [15] I. Cho, Graph }$W^{*}${\small -Probability Theory, (2004),
Preprint.}

{\small [16] I. Cho, Random Variables in Graph }$W^{*}${\small -Probability
Spaces over Diagonal Subalgebras, (2004), Preprint.}

{\small [17] P.\'{S}niady and R.Speicher, Continous Family of Invariant
Subspaces for R-diagonal Operators, Invent Math, 146, (2001) 329-363.}

{\small [18] R. Speicher, Combinatorial Theory of the Free Product with
Amalgamation and Operator-Valued Free Probability Theory, AMS Mem, Vol 132 ,
Num 627 , (1998).}

{\small [19] R. Speicher, Combinatorics of Free Probability Theory IHP
course note, available at www.mast.queensu.ca/\symbol{126}speicher.\strut }
\end{quote}

\end{document}